\documentclass[10pt]{amsart}
\usepackage{amsmath,amssymb,enumerate}
\usepackage{pstricks}

\newtheoremstyle{plainsl}%
	{\topsep}
	{\topsep}
	{\slshape} 
	{}
	{\normalfont\bfseries}
	{.}
	{ }
	{}
\swapnumbers

\newtheorem{theorem}{Theorem}[section]
\newtheorem{lemma}[theorem]{Lemma}
\newtheorem{corollary}[theorem]{Corollary}

\newcommand\cref[1]{Corollary~\ref{cor:#1}}
\newcommand\sref[1]{Section~\ref{sec:#1}}

\renewcommand\proof{\noindent\textsl{Proof. }}
\newcommand\sqr[2]{{\vbox{\hrule height.#2pt
    \hbox{\vrule width.#2pt height#1pt \kern#1pt
        \vrule width.#2pt}\hrule height.#2pt}}}
\renewcommand\qed{%
	\ifmmode\eqno\sqr53
	\else\nolinebreak\ \hfill\sqr53\medbreak\fi}

\newcommand{\HH}[2]{\ensuremath{H_{#1:#2}}}
\newcommand{\hh}{\ensuremath{\HH{n}{r}}}
\DeclareMathOperator{\dist}{dist}
\DeclareMathOperator{\diam}{diam}
\DeclareMathOperator{\aut}{Aut}

\begin{document}

\sloppy
\title{H\"{a}ggkvist--Hell Graphs: A class of Kneser-colorable graphs}
\author[Roberson]{David E. Roberson \\
                  \\
                  {D}\lowercase{epartment of} C\lowercase{ombinatorics and} O\lowercase{ptimization}\\
                  U\lowercase{niversity of} W\lowercase{aterloo} \\
                  200 U\lowercase{niversity} A\lowercase{venue} W\lowercase{est} \\
                  W\lowercase{aterloo}, ON N2L 3G1,
                  C\lowercase{anada} \\
                  \\
                  E\lowercase{mail}: \texttt{\lowercase{droberso@math.uwaterloo.ca}}}


\date{\today}

\begin{abstract}
For positive integers $n$ and $r$ we define the \emph{H\"{a}ggkvist--Hell graph}, \hh{}, to be the graph whose vertices are the ordered pairs $(h,T)$ where $T$ is an $r$-subset of $[n]$, and $h$ is an element of $[n]$ not in $T$. Vertices $(h_x,T_x)$ and $(h_y,T_y)$ are adjacent iff $h_x \in T_y$, $h_y \in T_x$, and $T_x \cap T_y = \varnothing$. These triangle-free arc transitive graphs are an extension of the idea of Kneser graphs, and there is a natural homomorphism from the H\"{a}ggkvist--Hell graph, \hh{}, to the corresponding Kneser graph, $K_{n:r}$. H\"{a}ggkvist and Hell introduced the $r=3$ case of these graphs, showing that a cubic graph admits a homomorphism to \HH{22}{3} if and only if it is triangle-free. Gallucio, Hell, and Ne\u{s}et\u{r}il also considered the $r=3$ case, proving that $\HH{n}{3}$ can have arbitrarily large chromatic number. In this paper we give the exact values for diameter, girth, and odd girth of all H\"{a}ggkvist--Hell graphs, and  we give bounds for independence, chromatic, and fractional chromatic number. Furthermore, we extend the result of Gallucio et al.~to any fixed $r \ge 2$, and we determine the full automorphism group of \hh{}, which is isomorphic to the symmetric group on $n$ elements.

\end{abstract}

\maketitle

\section{Introduction}
\label{sec:intro}

Let $n$ and $r$ be positive integers. Define the \emph{H\"{a}ggkvist--Hell graph}, \hh{}, to be the graph whose vertices are the ordered pairs $(h,T)$ where $T$ is an $r$-subset of $[n]$, and $h$ is an element of $[n]$ not in $T$. Vertices $(h_x,T_x)$ and $(h_y,T_y)$ are adjacent if $h_x \in T_y$, $h_y \in T_x$, and $T_x \cap T_y = \varnothing$. For a vertex $v = (h,T)$ of a H\"{a}ggkvist--Hell graph, we typically refer to $h$ as the \emph{head} of $v$ and $T$ as the \emph{tail} of $v$.

The reader may notice the similarity between \hh{} and the Kneser graph $K_{n:r}$, whose vertices are the $r$-subsets of $[n]$, and they are adjacent if disjoint. As we will see throughout this paper, this similarity is more than simply one of definition. In particular, the map that takes a vertex of \hh{} to its tail is a homomorphism to $K_{n:r}$. This homomorphism, and the relationship between \hh{} and $K_{n:r}$ in general, will serve as an important tool in the study of H\"{a}ggkvist--Hell graphs throughout this paper.

These graphs originally appeared in \cite{amote} in which H\"{a}ggkvist \& Hell showed that a cubic graph admits a homomorphism to \HH{22}{3} if and only if it is triangle-free. To our knowledge the only other reference to these graphs is \cite{color3} in which Galluccio, Hell, and Ne\u{s}et\u{r}il prove that \HH{n}{3} can have arbitrarily large chromatic number. Unlike the above two papers, which deal only with the $r=3$ case, we investigate properties of \hh{} in general.

%
%
%


In particular we give the diameter, girth and odd girth of \hh{} for all $n$ and $r$. We also give a lower bound on the size of independent sets in \hh{}. This bound is met with equality in all computed cases. We give two upper bounds on the chromatic number of \hh{}, one explicit and one recursive. We also show that for fixed $r \ge 2$, the chromatic number of \hh{} is unbounded, which extends the result of \cite{color3}. Using the bound on independent set size, we give an upper bound on the fractional chromatic number of \hh{} which shows that it is bounded for fixed $r$. Because of the relationship between fractional chromatic number and homomorphisms to Kneser graphs, the previous result implies the existence of a homomorphism from \hh{} to some Kneser graph which is not the homomorphism mentioned above. We then show how to construct this homomorphism using a technique that can be used to do the same for any vertex transitive graph. Lastly, we show that the automorphism group of \hh{} is isomorphic to the symmetric group on $n$ elements.

\section{Basic Properties}
\label{sec:basic}

Here we take time to remark on some easily seen properties of \hh{} before moving on to meatier fare. To begin, we note the following. 
\begin{itemize}
\item
\hh{} has $(r+1)\binom{n}{r+1} = (n-r)\binom{n}{r}$ vertices.
\item
\hh{} is regular with valency $r\binom{n-r-1}{r-1}$
\item
\hh{} is a subgraph of \HH{n'}{r} for $n \le n'$
\item
All H\"{a}ggkvist--Hell graphs are triangle-free.
\end{itemize}
The first two items follow from simple counting. The third is clear because \hh{} is simply the subgraph of \HH{n'}{r} induced by the vertices which only contain elements from $[n]$ in their head or tail. For the last item, note that all of the neighbors of a vertex $(h,T)$ have $h$ in their tail and therefore cannot be adjacent to each other.

As with the Kneser graphs, there are certain values of $n$ and $r$ for which \hh{} is not particularly interesting. We dispense of these cases now and deal with the more interesting cases in the rest of the paper. If $n < 2r$, then \hh{} can have no edges since there are no two disjoint sets of size $r$. Furthermore, if $n<r+1$, then \hh{} has no vertices. In the case when $r=1$, we see that every vertex has the form $(a,\{b\})$, which has only one neighbor, $(b,\{a\})$. So \HH{n}{1} is simply a matching. So throughout this paper we will assume that $r \ge 2$ and $n \ge 2r$.

In fact, if $n = 2r$ then \hh{} is the disjoint union of complete bipartite graphs. To see this note that any $r$-subset, $T$, of $[2r]$ is disjoint from exactly one other $r$-subset of $[2r]$, namely $\overline{T}$. Since the head of any vertex with tail $T$ must be an element of $\overline{T}$, and vice versa, every vertex with $T$ as its tail will be adjacent to exactly those vertices with $\overline{T}$ as their tail. Since there are $\frac{1}{2}\binom{2r}{r}$ pairs of such disjoint $r$-subsets, we have that
\[\HH{2r}{r} = \frac{1}{2}\binom{2r}{r}K_{r,r}.\]
This is similar to the case with Kneser graphs, which yields a matching of size $\frac{1}{2}\binom{2r}{r}$, except that we replace every edge with a $K_{r,r}$. This relationship between edges in $K_{n:r}$ and induced $K_{r,r}$'s in \hh{} holds for larger values of $n$, but in this case the vertices of two of these $K_{r,r}$'s may or may not intersect if the corresponding edges in $K_{n:r}$ intersect.

As with the Kneser graphs, the symmetric group on $n$ elements, $S_n$, is a subgroup of $\aut(\hh{})$. From the definition of adjacency, it is clear that $S_n$ acts arc transitively on \hh{}. In fact, $S_n$ is the entire automorphism group of \hh{} for $n \ge 2r+1$ which we will prove in \sref{aut}.

There are two important homomorphisms which arise quite naturally from the definition of \hh{}:
\begin{itemize}
\item
$(h,T) \mapsto h$ gives a homomorphism from \hh{} to the complete graph $K_n$.
\item
$(h,T) \mapsto T$ gives a homomorphism from \hh{} to the Kneser graph $K_{n:r}$.
\end{itemize}
The homomorphism to the Kneser graph is the one most used in this paper, though both homomorphisms relate to certain independent sets of \hh{}.

Since we know that \HH{n-1}{r} is a subgraph of \hh{}, one way to picture \hh{} is to think of \HH{n-1}{r} and then add the vertices with $n$ in their head or tail along with any necessary edges. In fact, it is helpful to distinguish between vertices with $n$ as their head, and vertices with $n$ in their tail. Doing this gives us the following partition which proves useful for studying many different aspects of H\"{a}ggkvist--Hell graphs.

\begin{figure}[h] 
   \centering
\psset{xunit=.5pt,yunit=.5pt,runit=.5pt}
\begin{pspicture}(574.18005371,193.37193298)
{
\newrgbcolor{curcolor}{1 1 1}
\pscustom[linestyle=none,fillstyle=solid,fillcolor=curcolor]
{
\newpath
\moveto(189.228725,94.58954541)
\curveto(189.228725,43.06611229)(158.06823405,1.24999266)(119.6740577,1.24999266)
\curveto(81.27988135,1.24999266)(50.1193904,43.06611229)(50.1193904,94.58954541)
\curveto(50.1193904,146.11297853)(81.27988135,187.92909816)(119.6740577,187.92909816)
\curveto(158.06823405,187.92909816)(189.228725,146.11297853)(189.228725,94.58954541)
\closepath
}
}
{
\newrgbcolor{curcolor}{0 0 0}
\pscustom[linewidth=2.49999994,linecolor=curcolor]
{
\newpath
\moveto(189.228725,94.58954541)
\curveto(189.228725,43.06611229)(158.06823405,1.24999266)(119.6740577,1.24999266)
\curveto(81.27988135,1.24999266)(50.1193904,43.06611229)(50.1193904,94.58954541)
\curveto(50.1193904,146.11297853)(81.27988135,187.92909816)(119.6740577,187.92909816)
\curveto(158.06823405,187.92909816)(189.228725,146.11297853)(189.228725,94.58954541)
\closepath
}
}
{
\newrgbcolor{curcolor}{1 1 1}
\pscustom[linestyle=none,fillstyle=solid,fillcolor=curcolor]
{
\newpath
\moveto(377.678822,98.52721851)
\curveto(377.678822,46.86293887)(356.24630304,4.93250901)(329.8383779,4.93250901)
\curveto(303.43045276,4.93250901)(281.9979338,46.86293887)(281.9979338,98.52721851)
\curveto(281.9979338,150.19149816)(303.43045276,192.12192801)(329.8383779,192.12192801)
\curveto(356.24630304,192.12192801)(377.678822,150.19149816)(377.678822,98.52721851)
\closepath
}
}
{
\newrgbcolor{curcolor}{0 0 0}
\pscustom[linewidth=2.5,linecolor=curcolor]
{
\newpath
\moveto(377.678822,98.52721851)
\curveto(377.678822,46.86293887)(356.24630304,4.93250901)(329.8383779,4.93250901)
\curveto(303.43045276,4.93250901)(281.9979338,46.86293887)(281.9979338,98.52721851)
\curveto(281.9979338,150.19149816)(303.43045276,192.12192801)(329.8383779,192.12192801)
\curveto(356.24630304,192.12192801)(377.678822,150.19149816)(377.678822,98.52721851)
\closepath
}
}
{
\newrgbcolor{curcolor}{1 1 1}
\pscustom[linestyle=none,fillstyle=solid,fillcolor=curcolor]
{
\newpath
\moveto(572.930062,96.10405851)
\curveto(572.930062,44.43977887)(551.49754304,2.50934901)(525.0896179,2.50934901)
\curveto(498.68169276,2.50934901)(477.2491738,44.43977887)(477.2491738,96.10405851)
\curveto(477.2491738,147.76833816)(498.68169276,189.69876801)(525.0896179,189.69876801)
\curveto(551.49754304,189.69876801)(572.930062,147.76833816)(572.930062,96.10405851)
\closepath
}
}
{
\newrgbcolor{curcolor}{0 0 0}
\pscustom[linewidth=2.5,linecolor=curcolor]
{
\newpath
\moveto(572.930062,96.10405851)
\curveto(572.930062,44.43977887)(551.49754304,2.50934901)(525.0896179,2.50934901)
\curveto(498.68169276,2.50934901)(477.2491738,44.43977887)(477.2491738,96.10405851)
\curveto(477.2491738,147.76833816)(498.68169276,189.69876801)(525.0896179,189.69876801)
\curveto(551.49754304,189.69876801)(572.930062,147.76833816)(572.930062,96.10405851)
\closepath
}
}
{
\newrgbcolor{curcolor}{0 0 0}
\pscustom[linewidth=2.5,linecolor=curcolor]
{
\newpath
\moveto(181.995906,141.18776298)
\curveto(181.995906,141.18776298)(221.762676,163.35411298)(255.055786,155.96533298)
\curveto(288.348896,148.57654298)(288.348896,148.57654298)(288.348896,148.57654298)
}
}
{
\newrgbcolor{curcolor}{0 0 0}
\pscustom[linestyle=none,fillstyle=solid,fillcolor=curcolor]
{
\newpath
\moveto(264.66096154,165.19630976)
\lineto(291.59962242,147.90021227)
\lineto(259.87332261,143.62372029)
\curveto(265.97273596,148.95194968)(267.87938711,157.67077713)(264.66096154,165.19630976)
\closepath
}
}
{
\newrgbcolor{curcolor}{0 0 0}
\pscustom[linewidth=2.5,linecolor=curcolor]
{
\newpath
\moveto(373.901706,131.83878298)
\curveto(373.901706,131.83878298)(413.927656,162.62385298)(447.437766,152.36216298)
\curveto(480.947866,142.10047298)(480.947866,142.10047298)(480.947866,142.10047298)
}
}
{
\newrgbcolor{curcolor}{0 0 0}
\pscustom[linestyle=none,fillstyle=solid,fillcolor=curcolor]
{
\newpath
\moveto(458.64112537,160.5326245)
\lineto(484.13529536,141.17045171)
\lineto(452.17086397,139.40364391)
\curveto(458.67062422,144.23541258)(461.25741578,152.77717938)(458.64112537,160.5326245)
\closepath
}
}
{
\newrgbcolor{curcolor}{0 0 0}
\pscustom[linewidth=2.5,linecolor=curcolor]
{
\newpath
\moveto(481.389046,58.79140298)
\curveto(481.389046,58.79140298)(441.189626,31.66674298)(409.015926,40.70830298)
\curveto(376.842236,49.74985298)(374.187636,50.74452298)(374.187636,50.74452298)
}
}
{
\newrgbcolor{curcolor}{0 0 0}
\pscustom[linestyle=none,fillstyle=solid,fillcolor=curcolor]
{
\newpath
\moveto(395.32240244,30.97948986)
\lineto(371.06321883,51.86820569)
\lineto(403.07584147,51.67205582)
\curveto(396.29208758,47.24785478)(393.18649419,38.88074295)(395.32240244,30.97948986)
\closepath
}
}
{
\newrgbcolor{curcolor}{0 0 0}
\pscustom[linewidth=2.5,linecolor=curcolor]
{
\newpath
\moveto(288.883536,56.83512298)
\curveto(288.883536,56.83512298)(248.140026,28.74644298)(215.167986,37.77814298)
\curveto(182.195936,46.80983298)(182.195936,46.80983298)(182.195936,46.80983298)
}
}
{
\newrgbcolor{curcolor}{0 0 0}
\pscustom[linestyle=none,fillstyle=solid,fillcolor=curcolor]
{
\newpath
\moveto(205.04193628,29.05046482)
\lineto(178.98221285,47.64447602)
\lineto(210.87981445,50.36283872)
\curveto(204.526896,45.33956301)(202.19574378,36.72451803)(205.04193628,29.05046482)
\closepath
}
}
{
\newrgbcolor{curcolor}{0 0 0}
\pscustom[linewidth=2.07423592,linecolor=curcolor]
{
\newpath
\moveto(54.355936,59.82652298)
\curveto(54.355936,59.82652298)(26.038363,48.17904298)(12.105483,64.22586298)
\curveto(3.564031,74.06325298)(-0.554124,88.99988298)(1.603056,102.36773298)
\curveto(2.964685,110.80561298)(16.53602,119.97336298)(26.038363,120.38971298)
\curveto(48.928096,121.39264298)(49.450326,121.39264298)(49.450326,121.39264298)
}
}
{
\newrgbcolor{curcolor}{0 0 0}
\pscustom[linestyle=none,fillstyle=solid,fillcolor=curcolor]
{
\newpath
\moveto(27.27581774,130.59625325)
\lineto(52.20494957,121.42917978)
\lineto(27.27581638,112.26210833)
\curveto(31.2584494,117.67434364)(31.23550138,125.07922087)(27.27581774,130.59625325)
\closepath
}
}

\uput[0](-40,150){$r\binom{n-r-2}{r-1}$}
\uput[0](180,180){$r\binom{n-r-2}{r-2}$}
\uput[0](385,180){$\binom{n-r-1}{r-1}$}
\uput[0](382,15){$r\binom{n-r-1}{r-1}$}
\uput[0](160,15){$(r-1)\binom{n-r-1}{r-1}$}
\uput[0](55,121){$\HH{n-1}{r}$}
\uput[0](294,121){$n \in T$}
\uput[0](490,121){$h = n$}
\uput[0](101,71){$C_1$}
\uput[0](309,71){$C_2$}
\uput[0](505,71){$C_3$}

\end{pspicture} 
   \caption{Diagram of 3-Cell Partition of \hh{}.}
   \label{fig:fig1}
\end{figure}
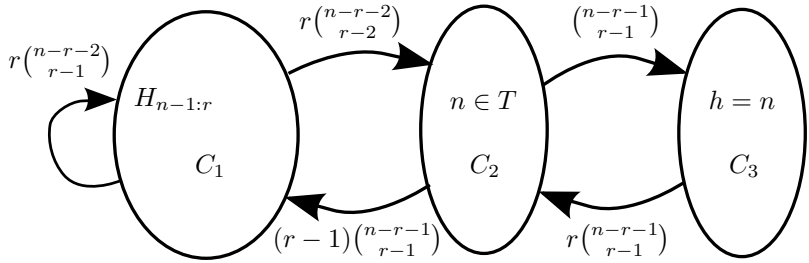

This partition turns out to be equitable, and the number on the arc from $C_i$ to $C_j$ is the number of neighbors a vertex in $C_i$ has in $C_j$. These can be determined by simple counting. Each cell contains a short description of the vertices it contains. So we see that $C_1$ is the \HH{n-1}{r} subgraph. The cell $C_2$ is the set of vertices with $n$ in their tail, which can also be described as the inverse image of a maximum independent set in $K_{n:r}$. Note that this set is independent in \hh{}. The cell $C_3$ is the set of vertices with $n$ as their head, which can also be described as the inverse image of a vertex (i.e. a maximum independent set) in $K_n$. Note that this set is independent in \hh{} and independent from the \HH{n-1}{r} subgraph.

Now that we have an idea of what \hh{} is like, we can start to ask some deeper questions about its structure. The first property of \hh{} we investigate is its diameter.

\section{Diameter}
\label{sec:diameter}

Let $G$ be a connected graph. The distance between two vertices $u$ and $v$ of $G$ is defined as the length of the shortest path in $G$ from $u$ to $v$, and is denoted $\dist(u,v)$. The diameter of $G$, $\diam(G)$, is defined as the maximum of $\dist(u,v)$ over all pairs of vertices, $u$ and $v$, in $G$.

The diameter of the Kneser graphs was given by Valencia-Pabon \& Vera in \cite{diameter}. This result will greatly aid us in determining the diameter of H\"{a}ggkvist--Hell graphs. The main technique is to use paths in $K_{n:r}$ to construct paths in \hh{}, and vice versa.

First we must determine the values of $n$ for which \hh{} is connected. For this we need the following lemma:

\begin{lemma}\label{tnected}
Any two vertices in \hh{} with the same tail are joined by a path.
\end{lemma}

\proof
Let $x = (h_x,T)$ and $y = (h_y,T)$, where $h_x \ne h_y$. Since $n \ge 2r$ and $h_x,h_y \notin T$, there exists an $r$-subset of $[n]\setminus T$ that contains both $h_x$ and $h_y$. Let $T'$ be such a set, and let $t \in T$. Then $x\sim(t,T')$ and $y\sim(t,T')$, and the path $x,(t',T),y$ connects $x$ and $y$.\qed

From this we can now prove that \hh{} is connected for $n \ge 2r+1$ using the analogous result for $K_{n:r}$. The following lemma gives the diameter of $K_{n:r}$ for $n \ge 2r+1$ which implies that it is connected. Though this is much more than we need here, we will need the full strength of this result later.

\begin{lemma}[Valencia-Pabon \& Vera]\label{knected}
For positive integers $n$ and $r$, $n \ge 2r+1$, the Kneser graph $K_{n:r}$ has diameter $\left\lceil \frac{r-1}{n-2r}\right\rceil + 1$.\qed
\end{lemma}

Combining these two results we get the following:

\begin{theorem}\label{connected}
For $n \ge 2r + 1$, the graph \hh{} is connected.
\end{theorem}

\proof
Let $x=(h_x,T_x)$ and $y=(h_y,T_y)$ be two vertices of \hh{}. Since $n \ge 2r + 1$, the graph $K_{n:r}$ is connected by Lemma~\ref{knected}, and so there is some path $T_x = T_0,T_1,\ldots,T_k = T_y$ in $K_{n:r}$. Now let $h_{i,1} \in T_{i-1}$ and $h_{i,2} \in T_{i+1}$ for all appropriate $i$. We see that $(h_{i,2},T_i)\sim (h_{i+1,1},T_{i+1})$ for all $i$ and by Lemma~\ref{tnected} $(h_{i,1},T_i)$ and $(h_{i,2},T_i)$ are joined by a path for all $i$, therefore $x$ and $y$ are joined by a path.\qed

Now that we know that \hh{} is connected for $n \ge 2r+1$, we can begin to speak of its diameter. Unlike Kneser graphs, which have diameter two for $n \ge 3r - 1$, H\"{a}ggkvist--Hell graphs have diameter at least four for all values of $n$.
\begin{lemma}\label{diaml4}
$\diam(\hh{}) \ge 4$.
\end{lemma}
\proof
Consider the vertices $x = (1,T_x)$ and $y = (1,T_y)$ such that $T_x \cap T_y = \varnothing$. Note that two such vertices always exist for $n \ge 2r+1$. We will show that $\dist(x,y) \ge 4$. Clearly $x$ and $y$ are not adjacent, since they have identical heads, thus they are not at distance one from each other. Now suppose that $x$ and $y$ share a common neighbor $z = (h_z,T_z)$. Then we have that $h_z \in T_x$ and $h_z \in T_y$, which is not possible since they are disjoint. Therefore $x$ and $y$ are at a distance of at least three from each other. Suppose that $\dist(x,y) = 3$. Then there exists two vertices $z_1 = (h_1,T_1)$ and $z_2 = (h_2,T_2)$ such that $P = x,z_1,z_2,y$ is a path. However, we see that this implies that $1 \in T_1$ and $1 \in T_2$, and therefore $T_1 \cap T_2 \ne \varnothing$ and so $z_1$ and $z_2$ are not adjacent and $P$ cannot be a path. Therefore $\dist(x,y) \ge 4$.\qed

For $n \ge \frac{5}{2}r$, this lower bound is achieved with equality.

\begin{lemma}\label{diam4}
For $n \ge \frac{5}{2}r$ the diameter of \hh{} is four.
\end{lemma}

\proof
Note that we only need to show that $\dist(x,y) \le 4$ for all vertices $x,y \in V(\hh{})$. Suppose that $n \ge \frac{5}{2}r$. Observe that this implies that $n \ge \left\lceil\frac{5}{2}r\right\rceil$, since $n$ is an integer. Let $x = (h_x,T_x)$ and $y = (h_y,T_y)$. We have four main cases:
\begin{enumerate}
\item
$h_x = h_y$;
\item
$h_x \in T_y$ and $h_y \notin T_x$;
\item
$h_x \in T_y$ and $h_y \in T_x$;
\item
$h_x \notin T_y$, $h_y \notin T_x$, and $h_x \ne h_y$.
\end{enumerate}
Let $C = T_x \cap T_y$ and $s = |C|$. Note that we must only consider $s \in \{0,1,\ldots,r-1\}$ for each of the above cases in order to prove our claim. If $s = r$, then either $x=y$ or $\dist(x,y) = 2$ by Lemma~\ref{tnected}, so we need not worry about these cases. We will prove the claim for the first case and then assure the reader that the other cases are similar.

Since $h_x = h_y$ in this case, we will refer to both as simply $h$. Let
\[D = [n]\setminus(T_x \cup T_y \cup \{h\}) = \{d_1,\ldots,d_{n-2r+s-1}\},\]
and $\ell = \left\lfloor\frac{r-s}{2}\right\rfloor$. Note that
\[T_x = \{x_1,\ldots,x_{r-s}\}\cup C \text{ and } T_y = \{y_1,\ldots,y_{r-s}\}\cup C\]
where $x_i \ne y_j$ for any $i,j$, and $r-s \ge 1$. Consider the vertices
\[z_x = (x_{r-s},\{h,y_1,\ldots,y_\ell,d_1,\ldots,d_{r-1-\ell}\})\]
\[z_y = (y_{r-s},\{h,x_1,\ldots,x_\ell,d_1,\ldots,d_{r-1-\ell}\})\]
\[w = (h,\{x_{r-s},y_{r-s},x_{\ell+1},\ldots,x_{r-s-1},y_{\ell+1},\ldots,y_{r-s-1}\}\cup C)\]
Note that the tail of $w$ has exactly $2r-s-2\ell = 2r-s-2\left\lfloor\frac{r-s}{2}\right\rfloor$ elements, which equals $r$ if $r-s$ is even and equals $r+1$ if $r-s$ is odd. But in the latter case we can just remove one of the elements from the tail that is not $x_{r-s}$ or $y_{r-s}$. Now it is straightforward to see that $P = x,z_x,w,z_y,y$ is a path from $x$ to $y$ as long as all of the indices are valid. Upon investigation one can see that the only thing we need to check is that $d_{r-1-\ell}$ exists, i.e. that $n-2r+s-1 \ge r-1-\ell$. However, this is equivalent to
\[n \ge 3r-s-\ell \ge 3r-s-\frac{r-s}{2} = \frac{5}{2}r-\frac{s}{2}.\]
Therefore, if $n \ge \frac{5}{2}r$, then $\dist(x,y) \le 4$. The other three cases are quite similar, so we will spare you the tedium.\qed

So we have determined the diameter of \hh{} for $n \ge \frac{5}{2}r$. In order to do the same for $2r+1 \le n < \frac{5}{2}r$ we must first prove the following two lower bounds.

\begin{lemma}\label{diamlb}
For $n \ge 2r+1$, the diameter of \hh{} is at least $\left\lceil \frac{r-1}{n-2r}\right\rceil + 1$.
\end{lemma}

\proof
Let $T_x,T_y \in V(K_{n:r})$ be such that $\dist(T_x,T_y) = \diam(K_{n:r})$. Let $P$ be a shortest path from $(1,T_x)$ to $(1,T_y)$ in \hh{}. Since the tails of consecutive vertices of $P$ must be disjoint, they represent a walk from $T_x$ to $T_y$ in $K_{n:r}$ of length equal to that of $P$. This implies that
\[\diam(K_{n:r}) = \dist_{K_{n:r}}(T_x,T_y) \le \dist_{\hh{}}((1,T_x),(1,T_y)) \le \diam(\hh{}).\] Since $\diam(K_{n:r}) = \left\lceil \frac{r-1}{n-2r}\right\rceil + 1$ by Lemma~\ref{knected}, the result is proven.\qed

\begin{lemma}\label{diamg4}
For $n < \frac{5}{2}r$, the diameter of \hh{} is strictly greater than four.
\end{lemma}

\proof
Consider the vertices $x = (h,T_x)$ and $y = (h,T_y)$ where $T_x \cap T_y = \varnothing$. Note such a pair of vertices exists for $n \ge 2r+1$. Let
\[T_x = \{x_1,\ldots,x_r\}, \ T_y = \{y_1,\ldots,y_r\} \text{, and } [n]\setminus (T_x \cup T_y \cup \{h\}) = D = \{d_1,\ldots,d_{k-1}\}\]
where $k = n-2r$. From the proof of Lemma~\ref{diaml4} we see that $\dist(x,y) \ge 4$, so we only need to show that there is no path of length four between $x$ and $y$. Suppose that $P = x,z_x,w,z_y,y$ is a path. We will show that we need at least $\frac{5}{2}r$ elements of $[n]$ for this path to exist. Immediately we see that $h \in T_{z_x}$ and $h \in T_{z_y}$, and WLOG we can say that $z_x = (x_r,T_{z_x})$ and $z_y = (y_r,T_{z_y})$.

We have two options for the head of $w$, either it is $h$, or it is some element of $D$. As it turns out, this does not make a difference, but for now we will assume that it is $h$. At the end of the proof we will show why the other case works out to be the same. Suppose $T_{z_x}$ and $T_{z_y}$ contain $i$ and $j$ elements from $D$ respectively. WLOG $i \le j$.

Suppose that $d \in T_{z_x} \cap D$ and $d \notin T_{z_y}$. Since $i \le j$, there must exist $d'\in T_{z_y} \cap D$ such that $d' \notin T_{z_x}$. Also, $d,d' \notin T_w$, but then we could simply replace the $d$ in $T_{z_x}$ with $d'$, and this will still be a path from $x$ to $y$ and it will use fewer elements from $[n]$, so we may assume that $(T_{z_x} \cap D) \subseteq T_{z_y}$. The other $r-i-1$ elements of $T_{z_x}$, and $r-j-1$ elements of $T_{z_y}$ come from $T_y\setminus y_r$ and $T_x\setminus x_r$ respectively.

So far, we have used $2r+1+j$ elements of $[n]$ in the vertices $x$, $y$, $z_x$, and $z_y$. Now we are left with the elements to be used in the tail of $w$. We know that $x_r,y_r \in T_w$, since these are the heads of $z_x$ and $z_y$ respectively. We are also able to use any of the other elements of $T_x \cup T_y$ not already used in $T_{z_x}$ or $T_{z_y}$, of which there are exactly
\[2r - 2 - (r-i-1) - (r-j-1) = i+j.\]
This leaves $r-2-i-j$ elements left in the tail of $w$, and these must come from $D\setminus T_{z_y}$. Thus we use a total of
\[(2r+1+j)+(r-2-i-j) = 3r-1-i\]
elements of $[n]$. However, this does not take into account the possibility that we were able to fill the tail of $w$ without using any elements of $D$, i.e. when $r-2-i-j \le 0$. In this case $3r-1-i \le 2r+1+j$, but we still use $2r+1+j$ elements of $[n]$ in our path $P$. In order to take this into account we must take the maximum of these two values. So the number of elements of $[n]$ that we use in the path $P$ is
\[\max\{3r-1-i,2r+1+j\}.\]
It is easy to see that letting $i=j$ can only reduce this maximum, and so we need to find the value of $i$ for which
\[\max\{3r-1-i,2r+1+i\}\]
is minimized. This will be minimized when $3r-1-i = 2r+1+i \Leftrightarrow i = \frac{1}{2}r - 1$. Note that $i$ must be an integer, but this can only increase the lower bound we get on $n$, and so we can ignore this. Plugging in this value of $i$ we see that we must use at least $\frac{5}{2}r$ elements of $[n]$ for the path $P$, thus proving the result. We see now that having an element of $D$ as the head of $w$ would have simply forced us to use an element of $T_{z_x} \cap D$, which only would have precluded us from having $i=0$, and would not have reduced the number of elements of $[n]$ that we needed for the path.\qed

The proof of the following lemma from \cite{diameter} is important for the proof of our final theorem on the diameter of \hh{}.

\begin{lemma}[Valencia-Pabon \& Vera]\label{dist}
Let $X,Y \in [n]^{(r)}$ be two different vertices in the Kneser graph $K_{n:r}$ with $2r+1\le n \le 3r-2$, such that $|X\cap Y| = s$. Then
\[\dist(X,Y) = \min\left\{2\left\lceil\frac{r-s}{n-2r}\right\rceil, 2\left\lceil\frac{s}{n-2r}\right\rceil + 1\right\}.\]
\end{lemma}

We only give the proof that this is an upper bound on $\dist(X,Y)$, because it is this portion that we will use for our proof of the diameter of the H\"{a}ggkvist--Hell graphs.
\vspace{1.7ex}

\proof
Let $k = n-2r$, so that $1 \le k < n-1$. Also, let $C = X\cap Y$, $s = |C|$, and $D = [n]\setminus(X\cup Y)$. Thus $|D|=s+k$.
Assume that $X=\{a_1,\ldots,a_{r-s}\}\cup C$, and $Y = \{b_1,\ldots,b_{r-s}\}\cup C$. Let $\ell = 2\left\lceil(r-s)/k\right\rceil$. Consider the path $X=T_0,T_1,\ldots,T_\ell = Y$ between $X$ and $Y$, where for $i < (r-s)/k$,
\begin{align*}
&T_{2i-1} = \{a_1,\ldots,a_{(i-1)k},b_{ik+1},\ldots,b_{r-s}\}\cup D, \\
&T_{2i} = \{b_1,\ldots,b_{ik},a_{ik+1},\ldots,a_{r-s}\}\cup C, \\
&\text{and} \\
&T_{\ell-1} = \{a_1,\ldots,a_{r-s-k}\}\cup D.
\end{align*}
Also, let $D' \subseteq D$ with $|D'| = s$. Consider the vertex $X' = (Y\setminus C)\cup D'$. Note that $X\cap X' = \varnothing$, and $s' = |X' \cap Y| = r-s$. Therefore, by the previous construction, there is a path between $X'$ and $Y$ with length equal to $2\left\lceil (r-s')/k \right\rceil = 2\left\lceil s/k \right\rceil$. Thus, there is a path between $X$ and $Y$ with length equal to $2\left\lceil s/k \right\rceil + 1$. So,
\[\dist(X,Y) \le \min\{2\left\lceil (r-s)/k \right\rceil, 2\left\lceil s/k \right\rceil + 1\}.\qed\]

We are now able to give the diameter of all connected H\"{a}ggkvist--Hell graphs in terms of their parameters.

\begin{theorem}\label{diameter}
For $n \ge \frac{5}{2}r$ the diameter of \hh{} is four. For $2r+1 \le n < \frac{5}{2}r$, the diameter of \hh{} is equal to $\max\left\{5,\left\lceil\frac{r-1}{n-2r}\right\rceil + 1\right\}$.
\end{theorem}

\proof
The first statement has already been proven as Lemma~\ref{diam4}. Also, Lemma~\ref{diamg4} and Lemma~\ref{diamlb} give the lower bound direction of the second statement. Thus we only have to show that we are able to achieve this bound for $n < \frac{5}{2}r$. We do this by showing that for any two vertices $(h_x,T_x)$ and $(h_y,T_y)$ in \hh{}, there is either a path between them of the same length as the shortest path between $T_x$ and $T_y$ in $K_{n:r}$, or there is a path between them of length at most 5. As in Lemma~\ref{diam4}, we have four main cases:
\begin{enumerate}
\item
$h_x = h_y$;
\item
$h_x \in T_y$ and $h_y \notin T_x$;
\item
$h_x \in T_y$ and $h_y \in T_x$;
\item
$h_x \notin T_y$, $h_y \notin T_x$, and $h_x \ne h_y$.
\end{enumerate}
We will use the same notation as in Lemma~\ref{dist}, so $k = n-2r$, $C = T_x \cap T_y$, $s=|C|$, $D = [n]\setminus(T_x\cup T_y)$, and $|D| = s+k$. Note that Lemma~\ref{tnected} takes care of the cases in which the vertices have the same tail. So we can assume that $s \le r-1$. We can also immediately take care of the cases with disjoint tails, as follows. Suppose that $T_x \cap T_y = \varnothing$. Then we have the following four cases:

\begin{itemize}
\item
If $h_x = h = h_y$, then let $t_x,t'_x \in T_x$ and $t_y,t'_y \in T_y$. Then there is a path of length five between $x$ and $y$ given by
\[(h,T_x),(t_x,\{h\}\cup (T_y \setminus t'_y)),(t_y,T_x),(t_x,T_y),(t_y,\{h\}\cup (T_x \setminus t'_x)),(h,T_y).\]

\item
If $h_x \in T_y$ and $h_y \notin T_x$, then let $t_x,t'_x \in T_x$ and $t_y \in T_y$. Then a path of length three from $x$ to $y$ is given by
\[(h_x,T_x),(t_x,T_y),(t_y,\{h_y\} \cup (T_x\setminus t'_x)),(h_y,T_y).\]

\item
If $h_x \in T_y$ and $h_y \in T_x$, then $x$ and $y$ are simply neighbors.
\vspace{1ex}

\item
If $h_x \notin T_y$, $h_y \notin T_x$, and $h_x \ne h_y$, then let $t_x,t'_x \in T_x$ and $t_y,t'_y \in T_y$. Then a path of length three between $x$ and $y$ is given by \[(h_x,T_x),(t_x,\{h_x\}\cup(T_y\setminus t'_y)),(t_y,\{h_y\} \cup (T_x\setminus t'_x)),(h_y,T_y).\]
\end{itemize}

Thus if the tails of two vertices are disjoint, then they are at a distance of at most five. Now we consider the cases where $r-k \le s \le r-1$. In this case we have that $|D| = s+k \ge r$.

\begin{enumerate}
\item
If $h_x = h = h_y$, then let $D' \subseteq D$ be such that $h \in D'$ and $|D'| = r$, and let $t \in T_x \cap T_y$. Then a path of length two between $x$ and $y$ is given by
\[(h,T_x),(t,D'),(h,T_y).\]

\item
If $h_x \in T_y$ and $h_y \notin T_x$, then let $D' \subseteq D$ be such that $h_y \in D'$ and $|D'| = r$, and let $d \in D'\setminus h_y$ and $t \in T_x \cap T_y$. Then there is a path of length four between $x$ and $y$ given by
\[(h_x,T_x),(t,\{h_x\}\cup(D'\setminus h_y)),(d,T_x),(t,D'),(h_y,T_y).\]

\item
If $h_x \in T_y$ and $h_y \in T_x$, let $D' \subseteq D$ be such that $|D'| = r-1$, and let $d \in D'$, $e \in D\setminus D'$, and $t \in T_x \cap T_y$. Then a path of length four between $x$ and $y$ is given by
\[(h_x,T_x),(t,\{h_x\}\cup D'),(d,\{e\}\cup(T_x \setminus h_y)),(t,\{h_y\}\cup D'),(h_y,T_y).\]

\item
If $h_x \notin T_y$, $h_y \notin T_x$, and $h_x \ne h_y$, then let $D' \subseteq D$ be such that $h_x,h_y \in D'$ and $|D'| = r$, and let $t \in T_x \cap T_y$. Then a path of length two between $x$ and $y$ is given by
\[(h,T_x),(t,D'),(h,T_y).\]
\end{enumerate}

So we have taken care of all cases in which $s \ge r-k$. For the remaining cases, we will be using the two paths between $T_x$ and $T_y$ given in the proof of Lemma~\ref{dist}. From them we construct two paths between $x$ and $y$ of lengths equal to those of the paths in the Kneser graph $K_{n:r}$. In order to do this we treat each vertex in the path of the Kneser graph as a tail of a vertex in \hh{}, and then we show that we are able to pick heads for each vertex in the interior of the path such that the adjacencies are preserved. After this, all that remains to show is that in each case we are able to choose the second and second to last vertices in the paths in the Kneser graph such that they contain the heads of $x$ and $y$ respectively.

For vertices in the interior of the paths this is trivial. Since $n < \frac{5}{2}r$, we have that if $T_1,T_2,T_3$ are three consecutive vertices in the path in $K_{n:r}$, then $T_1,T_3 \subseteq [n]\setminus T_2$ which has size less than $\frac{3}{2}r$, thus there must exist some element $t \in T_1 \cap T_3$, and we can pick this as the head of $T_2$ in \hh{}.

So all we need to show is that we are able to choose appropriate second and second to last vertices in the paths in the Kneser graph. We have to deal with each path separately:

Recall that $T_x = \{x_1,\ldots,x_{r-s}\}\cup C$ and $T_y = \{y_1,\ldots,y_{r-s}\}\cup C$. For the first path given in the proof of Lemma~\ref{dist}, we have that $T_1 = \{y_{k+1},\ldots,y_{r-s}\}\cup D$ and $T_{\ell-1} = \{x_1,\ldots,x_{r-s-k}\}\cup D$. Note that since $1 \le s \le r-k-1$, we have that $|D| = s+k \le r-1$. Now we go through the cases:

\vspace{1.7ex}

Case 1: $h_x = h_y$. In this case $h_x,h_y \in D \subseteq T_1,T_{\ell-1}$ and so we are done.

\vspace{1.7ex}

Case 2: $h_x \in T_y$ and $h_y \notin T_x$. Here $h_y \in D$, and if we let $y_{k+1} = h_x$, then $h_x \in T_1$ and $h_y \in T_{\ell-1}$.

\vspace{1.7ex}

Case 3: $h_x \in T_y$ and $h_y \in T_x$. Here if we let $y_{k+1} = h_x$ and $x_1 = h_y$, then $h_x \in T_1$ and $h_y \in T_{\ell-1}$.

\vspace{1.7ex}

Case 4: $h_x \notin T_y$, $h_y \notin T_x$, and $h_x \ne h_y$. Here we have that $h_x,h_y \in D \subseteq T_1,T_{\ell-1}$, and so we are done.

\vspace{1.7ex}

So we have shown that we can construct a path from $x$ to $y$ in \hh{} with the same length as the first path between $T_x$ and $T_y$ given in the proof of Lemma~\ref{dist}. Now we must do the same for the second path.

The second vertex in the second path is $T'_x = (T_y\setminus C)\cup D'$ where $D' \subseteq D$ such that $|D'| = s \ge 1$. Then, using the same construction as for the first path, let $C' = T'_x \cap T_y = T_y \setminus C$, and $E = [n]\setminus (T'_x \cup T_y) = [n]\setminus (T_y \cup D')$, and $|E| = r+k-s \ge 2k+1$. So $T'_x = \{d_1,\ldots,d_s\}\cup C'$ where $\{d_1,\ldots,d_s\} = D'$. Then the second to last vertex in the path is $T'_{\ell'-1} = \{d_1,\ldots,d_{s-k}\}\cup E$. Now we go through the cases:

\vspace{1.7ex}

Case 1: $h_x = h_y$. In this case $h_x,h_y \in D$, so if we let $h_x \in D'$ and $h_y \notin D'$ (possible since $|D\setminus D'| = k \ge 1$), then $h_x \in T'_x$ and $h_y \in E \subseteq T'_{\ell'-1}$ and so we are done.

\vspace{1.7ex}

Case 2: $h_x \in T_y$ and $h_y \notin T_x$. Here $h_x \in T_y \setminus C = C' \subseteq T'_x$, and we let $h_y \in D\setminus D' \subseteq E \subseteq T'_{\ell'-1}$. Thus we are done.

\vspace{1.7ex}

Case 3: $h_x \in T_y$ and $h_y \in T_x$. Here $h_x \in T_y \setminus C = C' \subseteq T'_x$, and $h_y \notin T_y$ and $h_y \notin D \supseteq D'$ and thus $h_y \in E \subseteq T'_{\ell'-1}$.

\vspace{1.7ex}

Case 4: $h_x \notin T_y$, $h_y \notin T_x$, and $h_x \ne h_y$. Here we have that $h_x,h_y \in D$ and we let $h_x \in D' \subseteq T'_x$ and $h_y \notin D'$. Thus $h_y \in E \subseteq T'_{\ell'-1}$.

\vspace{1.7ex}

Since one of these two paths must be a shortest path between $T_x$ and $T_y$ in $K_{n:r}$, we have now shown that for $n < \frac{5}{2}r$, any two vertices $x = (h_x,T_x)$ and $y = (h_y,T_y)$ of $\hh{}$ are either at a distance of at most five, or
\[\dist_{\hh{}}(x,y) \le \dist_{K_{n:r}}(T_x,T_y),\]
thus
\[\diam(\hh{}) \le \max\left\{5, \left\lceil\frac{r-1}{n-2r}\right\rceil + 1\right\}.\]
This completes the proof.\qed

This resolves all questions regarding the diameter of \hh{}.

\section{Odd Girth}
\label{sec:oddgirth}

The girth of a graph $G$ is defined as the length of the shortest cycle of $G$, whereas the odd girth of $G$ is likewise defined as the length of the shortest odd cycle of $G$. For H\"{a}ggkvist--Hell graphs the more interesting parameter turns out to be odd girth, since the girth is the same for any nonempty H\"{a}ggkvist--Hell graph. To see this, note that \hh{} contains \HH{2r}{r} as a subgraph, and this subgraph is a disjoint union of $K_{r,r}$'s and thus contains a four-cycle. Combine this with the fact that H\"{a}ggkvist--Hell graphs are triangle-free and we see that they always have girth four.

As with diameter, the odd girth of H\"{a}ggkvist--Hell graphs is closely related to the odd girth of Kneser graphs. However, also like diameter, equality does not always hold due to certain obstructions. In particular, the odd girth of \hh{} must be at least five since it is triangle-free. Fortunately, this seems to be the only obstruction. Similarly to diameter, the main technique we use in this section is to construct cycles in \hh{} using cycles in $K_{n:r}$ and vice versa. The following result by Poljak \& Tuza in \cite{odd_girth} gives the odd girth of the Kneser graphs which contain odd cycles.

\begin{theorem}[Poljak \& Tuza]\label{knog}
The odd girth of the Kneser graph $K_{n:r}$ is $2\left\lceil \frac{r}{n-2r} \right\rceil + 1$ for $n \ge 2r+1$.\qed
\end{theorem}

Since there is a homomorphism from \hh{} to $K_{n:r}$, the odd girth of \hh{} must be at least that of $K_{n:r}$. However, we give a direct proof as well.

\begin{lemma}\label{hhoglb}
For $n \ge 2r+1$, the odd girth of \hh{} is at least $2\left\lceil \frac{r}{n-2r} \right\rceil + 1$.
\end{lemma}

\proof
For $n \ge \frac{5}{2}r$, we have that $2\left\lceil \frac{r}{n-2r} \right\rceil + 1 \le 5$, which we have already established as a lower bound. So we can assume that $n \le \frac{5}{2}r$. Suppose that $C$ is a shortest odd cycle in \hh{}. If no tail is repeated in $C$, then the tails correspond to a cycle in $K_{n:r}$ and the result is proven by Theorem~\ref{knog}.

Otherwise, suppose that $x = (h_x,T)$ and $y = (h_y,T)$ are two vertices in $C$ with the same tail. Let $P$ be the path from $x$ to $y$ in $C$ with odd length. Now let $T_x$ be the tail of the unique neighbor of $x$ in $P$, and let $T_y$ be defined similarly. Since $n \le \frac{5}{2}r$, and $T_x,T_y \subseteq [n]\setminus T$, there must be an element $h$ of $[n]$ such that $h \in T_x \cap T_y$. Let $z = (h,T)$, and let $P'$ be the path $P$ with the ends, $x$ and $y$, removed. The cycle $C' = z,P',z$ is a shorter odd cycle than $C$, which is a contradiction.\qed

Now we are able to give the exact value of the odd girth of the H\"{a}ggkvist--Hell graphs for all values of $n$ and $r$.

\begin{theorem}\label{hhog}
For $n \ge 2r+1$, the odd girth of \hh{} is
\[\max\left\{5, 2\left\lceil \frac{r}{n-2r} \right\rceil + 1\right\}.\]
\end{theorem}

\proof
Lemma~\ref{hhoglb} and the above gives the lower bound direction, so we only need to show that it can be achieved. First we consider the case where $n < 3r$. In this case $2\left\lceil \frac{r}{n-2r} \right\rceil + 1 \ge 5$, so we will show that we can obtain an odd cycle with this length.

Consider a shortest odd cycle $C$ in $K_{n:r}$, this has length $2\left\lceil \frac{r}{n-2r} \right\rceil + 1$. We will view the vertices of $C$ as tails and show that we can pick a head for each so that the adjacencies in $C$ are preserved. Consider a vertex $T$ in $C$, with neighbors $T_1$ and $T_2$ in $C$. Since $n < 3r$, and $T_1,T_2 \subseteq [n]\setminus T$, we have that $T_1$ and $T_2$ cannot be disjoint. So let $h \in T_1 \cap T_2$ and let this be the head of $T$. If we do this for each tail then we will have a cycle $C'$ in \hh{} with the same length as $C$.

Now we still need to consider the case where $n \ge 3r$. However, for $\frac{5}{2}r \le n \le 3r-1$, we have that $2\left\lceil \frac{r}{n-2r} \right\rceil + 1 = 5$, and so for $n = 3r-1$, the odd girth of \hh{} is five, which is as small as possible. Now for $n \ge 3r$, \hh{} contains \HH{3r-1}{r} as a subgraph, which means that it has odd girth at most five. But then it must have odd girth exactly five, and $2\left\lceil \frac{r}{n-2r} \right\rceil + 1 = 3$ for $n \ge 3r$, which proves the result.\qed

Note that we actually need to be somewhat careful in the above proof when saying that \HH{3r-1}{r} has odd girth five, since we need that $3r-1 \ge \frac{5}{2}r$. But this is in fact always true for $r \ge 2$.

\section{Subgraphs}
\label{sec:subgraphs}

The proofs above have shown us that we are sometimes able to use subgraphs in $K_{n:r}$ to construct isomorphic copies of these subgraphs in \hh{}. The next theorem considers this a bit more generally.

\begin{theorem}\label{subgraphs}
Let $n \ge 2r+1$. For any subgraph $G$ of $K_{n:r}$ with maximum degree strictly less than $\frac{n-r}{n-2r}$, there is a subgraph of \hh{} isomorphic to $G$.
\end{theorem}

\proof
Let $G$ be a subgraph of $K_{n:r}$, let $\Delta$ be the maximum degree of $G$, and suppose that $\Delta < \frac{n-r}{n-2r}$. Consider a vertex $X$ of $G$. Since $X$ is a vertex of $K_{n:r}$, it is an $r$-subset of $[n]$ whose neighbors in $G$ are also $r$-subsets of $[n]$, and they are disjoint from $X$. It suffices to show that we can pick a head $h_X$ for $X$ such that $h_X$ is an element of every neighbor of $X$ in $G$. Doing this for all vertices in $G$ completes the proof. So we must show that the neighbors of $X$ all share a common element, then we can pick that element and we are done. Let $k$ be the degree of $X$ in $G$. Now since the neighbors of $X$ are all disjoint from $X$, they must all draw their elements from the same set of size $n-r$. Let us call this set $S$. Now suppose that no element of $S$ is common to all of the neighbors of $X$. Then each element of $S$ is in at most $k - 1$ neighbors of $X$. Viewing this as a $r$-uniform hypergraph with $S$ as the ground set and the $k$ neighbors of $X$ as the hyperedges, we know that the degree sum of the vertices (at most $(k-1)(n-r)$) is equal to the degree sum of the hyperedges ($kr$), thus we have the following string of inequalities:
\begin{align*}
(k-1)(n-r) &\ge kr \\
kn-kr-n+r &\ge kr \\
kn-2kr &\ge n-r \\
k(n-2r) &\ge n-r
\end{align*}
Therefore, 
\[k \ge \frac{n-r}{n-2r} > \Delta\]
which is a contradiction. Therefore the neighbors of $X$ must share a common element, and we are done.\qed

Note that in the case of the odd graphs, when $n = 2r+1$, the condition we get is $\Delta < r+1$, i.e. $\Delta \le r$. However, in this case the valency of $K_{2r+1:r}$ is $r+1$, and since Kneser graphs are vertex transitive, they either have a perfect matching or a matching missing exactly one vertex. Therefore, \hh{} contains a copy of $K_{2r+1:r}$ minus a perfect matching whenever $\binom{2r+1}{r}$ is even, and contains a copy of $K_{2r+1:r}$ minus a maximum matching and one edge incident to the single vertex the matching misses whenever $\binom{2r+1}{r}$ is odd.

\section{Independent Sets}
\label{sec:indysets}

In this section we find large independent sets of \hh{} and conjecture as to their maximality. For small $n$ it seems that the best we can do is take inverse images of maximum independent sets of $K_{n:r}$, i.e. sets of vertices with a common element in their tail. However, as $n$ increases, the set of vertices whose head is larger than any element in their tail outgrows this set and is also independent. There is a more general formulation of this second set related to inverse images of vertices (i.e. independent sets) in $K_n$. With these two ideas in hand, we figure out the optimal way of combining them to produce a larger independent set. This gives us a lower bound on the independence number of \hh{} which is met with equality in all computed cases. We use $\alpha(X)$ to denote the size of a largest independent set of the graph $X$.

Our main tool for finding independent sets will be the 3-cell partition of \hh{} described in \sref{basic}. We give the diagram again for the reader's ease of reference.
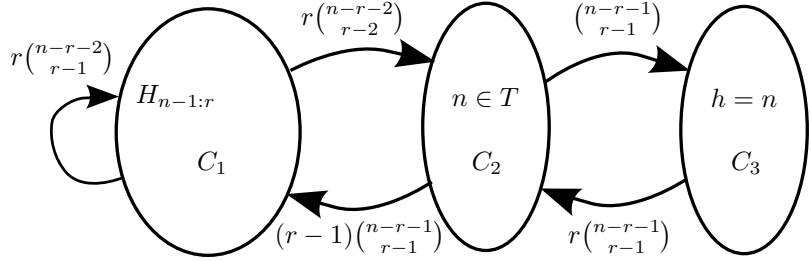
\begin{figure}[h] 
   \centering
\psset{xunit=.5pt,yunit=.5pt,runit=.5pt}
\begin{pspicture}(574.18005371,193.37193298)
{
\newrgbcolor{curcolor}{1 1 1}
\pscustom[linestyle=none,fillstyle=solid,fillcolor=curcolor]
{
\newpath
\moveto(189.228725,94.58954541)
\curveto(189.228725,43.06611229)(158.06823405,1.24999266)(119.6740577,1.24999266)
\curveto(81.27988135,1.24999266)(50.1193904,43.06611229)(50.1193904,94.58954541)
\curveto(50.1193904,146.11297853)(81.27988135,187.92909816)(119.6740577,187.92909816)
\curveto(158.06823405,187.92909816)(189.228725,146.11297853)(189.228725,94.58954541)
\closepath
}
}
{
\newrgbcolor{curcolor}{0 0 0}
\pscustom[linewidth=2.49999994,linecolor=curcolor]
{
\newpath
\moveto(189.228725,94.58954541)
\curveto(189.228725,43.06611229)(158.06823405,1.24999266)(119.6740577,1.24999266)
\curveto(81.27988135,1.24999266)(50.1193904,43.06611229)(50.1193904,94.58954541)
\curveto(50.1193904,146.11297853)(81.27988135,187.92909816)(119.6740577,187.92909816)
\curveto(158.06823405,187.92909816)(189.228725,146.11297853)(189.228725,94.58954541)
\closepath
}
}
{
\newrgbcolor{curcolor}{1 1 1}
\pscustom[linestyle=none,fillstyle=solid,fillcolor=curcolor]
{
\newpath
\moveto(377.678822,98.52721851)
\curveto(377.678822,46.86293887)(356.24630304,4.93250901)(329.8383779,4.93250901)
\curveto(303.43045276,4.93250901)(281.9979338,46.86293887)(281.9979338,98.52721851)
\curveto(281.9979338,150.19149816)(303.43045276,192.12192801)(329.8383779,192.12192801)
\curveto(356.24630304,192.12192801)(377.678822,150.19149816)(377.678822,98.52721851)
\closepath
}
}
{
\newrgbcolor{curcolor}{0 0 0}
\pscustom[linewidth=2.5,linecolor=curcolor]
{
\newpath
\moveto(377.678822,98.52721851)
\curveto(377.678822,46.86293887)(356.24630304,4.93250901)(329.8383779,4.93250901)
\curveto(303.43045276,4.93250901)(281.9979338,46.86293887)(281.9979338,98.52721851)
\curveto(281.9979338,150.19149816)(303.43045276,192.12192801)(329.8383779,192.12192801)
\curveto(356.24630304,192.12192801)(377.678822,150.19149816)(377.678822,98.52721851)
\closepath
}
}
{
\newrgbcolor{curcolor}{1 1 1}
\pscustom[linestyle=none,fillstyle=solid,fillcolor=curcolor]
{
\newpath
\moveto(572.930062,96.10405851)
\curveto(572.930062,44.43977887)(551.49754304,2.50934901)(525.0896179,2.50934901)
\curveto(498.68169276,2.50934901)(477.2491738,44.43977887)(477.2491738,96.10405851)
\curveto(477.2491738,147.76833816)(498.68169276,189.69876801)(525.0896179,189.69876801)
\curveto(551.49754304,189.69876801)(572.930062,147.76833816)(572.930062,96.10405851)
\closepath
}
}
{
\newrgbcolor{curcolor}{0 0 0}
\pscustom[linewidth=2.5,linecolor=curcolor]
{
\newpath
\moveto(572.930062,96.10405851)
\curveto(572.930062,44.43977887)(551.49754304,2.50934901)(525.0896179,2.50934901)
\curveto(498.68169276,2.50934901)(477.2491738,44.43977887)(477.2491738,96.10405851)
\curveto(477.2491738,147.76833816)(498.68169276,189.69876801)(525.0896179,189.69876801)
\curveto(551.49754304,189.69876801)(572.930062,147.76833816)(572.930062,96.10405851)
\closepath
}
}
{
\newrgbcolor{curcolor}{0 0 0}
\pscustom[linewidth=2.5,linecolor=curcolor]
{
\newpath
\moveto(181.995906,141.18776298)
\curveto(181.995906,141.18776298)(221.762676,163.35411298)(255.055786,155.96533298)
\curveto(288.348896,148.57654298)(288.348896,148.57654298)(288.348896,148.57654298)
}
}
{
\newrgbcolor{curcolor}{0 0 0}
\pscustom[linestyle=none,fillstyle=solid,fillcolor=curcolor]
{
\newpath
\moveto(264.66096154,165.19630976)
\lineto(291.59962242,147.90021227)
\lineto(259.87332261,143.62372029)
\curveto(265.97273596,148.95194968)(267.87938711,157.67077713)(264.66096154,165.19630976)
\closepath
}
}
{
\newrgbcolor{curcolor}{0 0 0}
\pscustom[linewidth=2.5,linecolor=curcolor]
{
\newpath
\moveto(373.901706,131.83878298)
\curveto(373.901706,131.83878298)(413.927656,162.62385298)(447.437766,152.36216298)
\curveto(480.947866,142.10047298)(480.947866,142.10047298)(480.947866,142.10047298)
}
}
{
\newrgbcolor{curcolor}{0 0 0}
\pscustom[linestyle=none,fillstyle=solid,fillcolor=curcolor]
{
\newpath
\moveto(458.64112537,160.5326245)
\lineto(484.13529536,141.17045171)
\lineto(452.17086397,139.40364391)
\curveto(458.67062422,144.23541258)(461.25741578,152.77717938)(458.64112537,160.5326245)
\closepath
}
}
{
\newrgbcolor{curcolor}{0 0 0}
\pscustom[linewidth=2.5,linecolor=curcolor]
{
\newpath
\moveto(481.389046,58.79140298)
\curveto(481.389046,58.79140298)(441.189626,31.66674298)(409.015926,40.70830298)
\curveto(376.842236,49.74985298)(374.187636,50.74452298)(374.187636,50.74452298)
}
}
{
\newrgbcolor{curcolor}{0 0 0}
\pscustom[linestyle=none,fillstyle=solid,fillcolor=curcolor]
{
\newpath
\moveto(395.32240244,30.97948986)
\lineto(371.06321883,51.86820569)
\lineto(403.07584147,51.67205582)
\curveto(396.29208758,47.24785478)(393.18649419,38.88074295)(395.32240244,30.97948986)
\closepath
}
}
{
\newrgbcolor{curcolor}{0 0 0}
\pscustom[linewidth=2.5,linecolor=curcolor]
{
\newpath
\moveto(288.883536,56.83512298)
\curveto(288.883536,56.83512298)(248.140026,28.74644298)(215.167986,37.77814298)
\curveto(182.195936,46.80983298)(182.195936,46.80983298)(182.195936,46.80983298)
}
}
{
\newrgbcolor{curcolor}{0 0 0}
\pscustom[linestyle=none,fillstyle=solid,fillcolor=curcolor]
{
\newpath
\moveto(205.04193628,29.05046482)
\lineto(178.98221285,47.64447602)
\lineto(210.87981445,50.36283872)
\curveto(204.526896,45.33956301)(202.19574378,36.72451803)(205.04193628,29.05046482)
\closepath
}
}
{
\newrgbcolor{curcolor}{0 0 0}
\pscustom[linewidth=2.07423592,linecolor=curcolor]
{
\newpath
\moveto(54.355936,59.82652298)
\curveto(54.355936,59.82652298)(26.038363,48.17904298)(12.105483,64.22586298)
\curveto(3.564031,74.06325298)(-0.554124,88.99988298)(1.603056,102.36773298)
\curveto(2.964685,110.80561298)(16.53602,119.97336298)(26.038363,120.38971298)
\curveto(48.928096,121.39264298)(49.450326,121.39264298)(49.450326,121.39264298)
}
}
{
\newrgbcolor{curcolor}{0 0 0}
\pscustom[linestyle=none,fillstyle=solid,fillcolor=curcolor]
{
\newpath
\moveto(27.27581774,130.59625325)
\lineto(52.20494957,121.42917978)
\lineto(27.27581638,112.26210833)
\curveto(31.2584494,117.67434364)(31.23550138,125.07922087)(27.27581774,130.59625325)
\closepath
}
}

\uput[0](-40,150){$r\binom{n-r-2}{r-1}$}
\uput[0](180,180){$r\binom{n-r-2}{r-2}$}
\uput[0](385,180){$\binom{n-r-1}{r-1}$}
\uput[0](382,15){$r\binom{n-r-1}{r-1}$}
\uput[0](160,15){$(r-1)\binom{n-r-1}{r-1}$}
\uput[0](55,121){$\HH{n-1}{r}$}
\uput[0](294,121){$n \in T$}
\uput[0](490,121){$h = n$}
\uput[0](101,71){$C_1$}
\uput[0](309,71){$C_2$}
\uput[0](505,71){$C_3$}

\end{pspicture} 
   \caption{Diagram of 3-Cell Partition of \hh{}.}
   \label{fig:fig1}
\end{figure}

The following two results are clear from the above diagram.

\begin{theorem}\label{indepkneserbd}
For $n \ge 2r$,
\[\alpha(\hh) \ge r\binom{n-1}{r} = (n-r)\binom{n-1}{r-1}.\]
\end{theorem}

\proof
An independent set of this size is given by the inverse image of a maximum independent set in the Kneser graph $K_{n:r}$. Equivalently, this is the set of all vertices of \hh{} with a common element in their tail i.e.~the middle cell in the above partition.\qed

\begin{theorem}\label{indeprecursivebd}
For all $r \ge 2$,
\[\alpha(\hh) \ge \alpha(\HH{n-1}{r}) + \binom{n-1}{r}.\]
\end{theorem}

\proof
The set of vertices of \hh{} that have $n$ as their head (the rightmost cell above) is an independent set of size $\binom{n-1}{r}$. This set is also independent from the \HH{n-1}{r} subgraph formed by vertices of \hh{} that do not contain $n$ in their head or tail (the leftmost cell above).\qed

The last theorem immediately gives us the following corollary:

\begin{corollary}
\[\alpha(\hh) \ge \sum_{i=r+1}^{n}\binom{i-1}{r} = \binom{n}{r+1}\].
\end{corollary}

\proof
An independent set of this size is obtained by recursively using the bound above, i.e. taking the vertices of \HH{k}{r} with $k$ as their head for $k = n,n-1,\ldots,1$. This independent set clearly has size equal to the summation above; the following equality is known but also follows easily from the realization that this independent set is exactly the set of vertices of \hh{} whose head is greater than any element in their tail. Since every subset of $r+1$ elements of $[n]$ gives rise to exactly one vertex of \hh{} in this set, it must have size $\binom{n}{r+1}$.\qed

Note that the above two lower bounds on the independence number of \hh{} can be rewritten as $\frac{r}{n}|V(\hh)|$ and $\frac{1}{r+1}|V(\hh)|$ respectively. We will frequently refer to these two types of independent sets as Kneser-type and recursive-type independent sets.

From the above we see that for $n \le r^2 + r$, the Kneser type independent set is larger than the recursive type independent set. This turns out to be to our advantage, since it allows us to start constructing a large independent set using the recursive approach, then stop and take the Kneser type independent set of the remaining H\"{a}ggkvist--Hell graph to obtain an independent set larger than either of our more pure-bred types. The key, then, is determining the optimal point to switch from one approach to the other. We proceed to this end.

We begin with the usual, requisite notation. For fixed $r$ we define the following:
\begin{align*}
\mathcal{H}_i(j) &= \{x \in V(\HH{i}{r}) : h_x = j\} \\ 
\mathcal{T}_i(j) &= \{x \in V(\HH{i}{r}) : T_x \ni j\} 
\end{align*}
In other words, $\mathcal{H}_i(j)$ is the part of a recursive type independent set we get from one recursive step; while $\mathcal{T}_i(j)$ is a Kneser type independent set of \HH{i}{r}. These definitions allow us to more concisely describe the independent sets we will be constructing. We can now formally define $\alpha'(\hh{})$, which is informally the size of the largest independent set we can construct by combining our above two ideas.
\begin{align*}
\alpha'(\hh{}) &= \max\{|\mathcal{T}_n(n)|, |\mathcal{H}_n(n)| + \alpha'(\HH{n-1}{r})\} \\
&= \max\left\{r\binom{n-1}{r},\binom{n-1}{r} + \alpha'(\HH{n-1}{r})\right\}
\end{align*}

Note that we take $\alpha'(\HH{2r}{r})$ to be $|\mathcal{T}_{2r}(2r)| = \frac{1}{2}|V(\HH{2r}{r})|$, which is maximum since $\HH{2r}{r}$ is the disjoint union of complete bipartite graphs. The next theorem states that the best point to stop using the recursive technique and take the Kneser type independent set of the remaining graph is at either of $r^2$ and $r^2 + 1$.

\begin{theorem}
\label{bestindybd}
For $2r \le n \le r^2 + 1$, $\alpha'(\hh{}) = r\binom{n-1}{r}$ and $\mathcal{T}_n(1)$ is an independent set of this size. For $n \ge r^2$, $\alpha'(\hh{}) = \binom{n}{r+1} + \frac{r-1}{r+1}\binom{r^2}{r}$ and
\[\left(\bigcup_{i=r^2+1}^{n} \mathcal{H}_i(i)\right) \cup \mathcal{T}_{r^2}(1)\]
is an independent set of this size.
\end{theorem}

\proof
We first determine, for a given $n$, whether it is better to take the Kneser type independent set or recurse once and then take the Kneser type independent set of the remaining \HH{n-1}{r} subgraph. This turns out to be the only case we need to consider explicitly. The size of a Kneser type independent set of \hh{} is $r\binom{n-1}{r}$. The set obtained by recursing once and then taking a Kneser type set is
\[\binom{n-1}{r} + r\binom{n-2}{r} = \left(1+r\frac{n-r-1}{n-1}\right)\binom{n-1}{r}\]
So we just need to compare $r$ and $1+r\frac{n-r-1}{n-1}$. We consider their ratio:
\begin{align*}
\frac{1+r\frac{n-r-1}{n-1}}{r} &= \frac{n-1 + nr - r^2 - r}{nr-r} \\
&= 1 + \frac{n-(r^2 + 1)}{r(n-1)}
\end{align*}
So we see that it is better to recurse once and then take the Kneser type independent set if and only if $n > r^2 + 1$. But of course this means that if we have $n > r^2 + 1$ then we should continue to recurse at least until we have a \HH{r^2 + 1}{r} remaining.

From the above we know that it is better to take a Kneser type independent set from \HH{r^2}{r} rather than recurse once more and then stop. But perhaps it is even better to recurse more than once and then take a Kneser type independent set in the remaining \HH{k}{r}, where $k < r^2 - 1$. If we choose $k$ to be as large as possible for this to occur, then it must have been better to recurse once at $k + 1$ than stopping there. However this implies that $k + 1 > r^2 + 1$ which is not the case.

As we can see from above, for \HH{r^2 + 1}{r} both recursing once and not recursing at all give independent sets of the same size. This explains the overlap in the conditions on $n$ in the theorem statement. So we see that for $n \ge r^2 + 1$, the best strategy is to stop recursing at either of $r^2$ and $r^2 + 1$. The size of this independent set will then be the size of a recursive type independent set of \hh{} plus the difference between the size of a Kneser type independent set and recursive type independent set of $\HH{r^2 + 1}{r}$:
\begin{align*}
\binom{n}{r+1} + r\binom{r^2 + 1 - 1}{r} - \binom{r^2 + 1}{r+1} &= \binom{n}{r+1} + \left(r - \frac{r^2 + 1}{r+1}\right)\binom{r^2}{r} \\
&= \binom{n}{r+1} + \frac{r-1}{r+1}\binom{r^2}{r}
\end{align*}
For $n \le r^2$ the best strategy is to simply take a Kneser type independent set which has size $r\binom{n-1}{r}$.\qed

We would like to be able to say that the independent sets from the above theorem are maximum independent sets, but presently a proof of this still eludes us. However we are at least able to say that they are maximal.

\begin{theorem}\label{bigindepmaximal}
The sets $\mathcal{T}_n(1)$ and $\left(\bigcup_{i=r^2+1}^{n} \mathcal{H}_i(i)\right) \cup \mathcal{T}_{r^2}(1)$ are maximal.
\end{theorem}

\proof
We give only an outline, but it is not difficult. $\mathcal{T}_i(1)$ is maximal for all $i \ge 2r$. The union of $\mathcal{H}_n(n)$ and a maximal independent set of \HH{n-1}{r} is maximal. This is clear from the 3-cell partition.\qed

Note that if $n \ge r^2 + 1$, then the independent set we construct with this approach has size $\frac{|V(\HH{n}{r})|}{r+1} + \frac{r-1}{r+1}\binom{r^2}{r}$. This is an improvement of $\frac{r-1}{r+1}\binom{r^2}{r}$ over the recursive type independent sets, which are larger than the Kneser type independent sets for these values of $n$.
Somewhat surprisingly, for $n \ge r^2 + 1$ it is possible to find two independents set of this size which are disjoint.

\begin{theorem}\label{twobigsets}
For $n \ge r^2 + 1$, there exist two disjoint independent sets in $\HH{n}{r}$, both having size $\binom{n}{r+1} + \frac{r-1}{r+1}\binom{r^2}{r}$.
\end{theorem}

\proof
Define $\sigma \in S_n$ to be the permutation define by $\sigma(i) = n+1-i$. The first independent set is the same as in the theorem above, $\left(\bigcup_{i=r^2+1}^{n} \mathcal{H}_i(i)\right) \cup \mathcal{T}_{r^2}(1)$, we will refer to this set as $\mathcal{L}^+$. The other independent set of this size is $\mathcal{L}^- = \sigma(\mathcal{L}^+)$. These two sets clearly have the correct size, all that is left is to show that they are disjoint.

Note that $\mathcal{L}^+$ is the set of vertices with heads from $\{r^2+1,\ldots,n\}$ that are greater than any element in their tails, call this $\mathcal{M}^+$, along with the vertices using only elements from $\{1,\ldots,r^2\}$ with $1$ in their tails, call this $\mathcal{N}^+$. Similarly, $\mathcal{L}^-$ is the set of vertices with heads from $\{1,\ldots,n-r^2\}$ that are less than any element in their tails, call this $\mathcal{M}^-$, along with the vertices using only elements from $\{n-r^2+1,\ldots,n\}$ with $n$ in their tails, call this $\mathcal{N}^-$. Clearly $\mathcal{M}^+$ and $\mathcal{M}^-$ are disjoint since the head of a vertex cannot be both larger and smaller than everything in its tail. $\mathcal{M}^+$ and $\mathcal{N}^-$ are disjoint because if $n$ is in the tail of a vertex, then its head cannot be larger than everything in its tail. Similarly, $\mathcal{M}^-$ and $\mathcal{N}^+$ are disjoint. Now all that is left is to show that $\mathcal{N}^+$ and $\mathcal{N}^-$ are disjoint. But the vertices in $\mathcal{N}^+$ only use elements from $\{1,\ldots,r^2\}$, and $n \ge r^2+1$, and so they cannot have $n$ in their tails. Therefore, $\mathcal{M}^+ \sqcup \mathcal{N}^+ = \mathcal{L}^+$ and $\mathcal{M}^- \sqcup \mathcal{N}^- = \mathcal{L}^-$ are disjoint.\qed

Note that we can always obtain two disjoint independent sets of size $\frac{|V(\hh{})|}{r+1}$: the sets of vertices whose head is greater/less than any element is its tail. This means we can color $\frac{2}{r+1}$ of the vertices using only two colors. We will investigate the chromatic number of \hh{} more in the next chapter.

We conclude our study of the independent sets of \hh{} with some values of $\alpha(\hh{})$ which we have computed for small values of $r$ and $n$.

\begin{center}
	\begin{tabular}{|c|c|c|c|}
    \hline $r$ & $n$ & $\alpha(H_{n:r})$ & $|V(H_{n:r})|$ \\
    \hline
		2 & 4 & 6 & 12 \\
		2 & 5 & 12 & 30 \\
		2 & 6 & 22 & 60 \\
		2 & 7 & 37 & 105 \\
		2 & 8 & 58 & 168 \\
		\hline
		3 & 6 & 30 & 60 \\
		3 & 7 & 60 & 140 \\
		3 & 8 & 105 & 280 \\
		\hline
	\end{tabular} \\
	\textsc{Table 1.} Independence numbers of \hh{}.
\end{center}
\addtocounter{table}{1}

Note that all values of $\alpha(\hh{})$ that we have computed agree with our best lower bound. In \sref{frac} we will use some of these results to give us bounds on the fractional chromatic number of \hh{}.

\section{Chromatic Number}
\label{sec:chromnum}

The chromatic number of the Kneser graph, $K_{n:r}$, is $n-2r+2$ for $n \ge 2r$. This easy upper bound was shown to hold with equality by Lov\'{a}sz \cite{lovasz}. We show that this is also an upper bound for the chromatic number of \hh{}, but it is not always met with equality. In addition to this, we give a recursive bound on $\chi(\hh{})$, showing that it increases by at most one when $n$ is increased by one. We also compute $\chi(\hh{})$ for small $n$ and $r$, which is how we conclude that the Kneser bound is not always met with equality. Galluccio et al.~\cite{color3} show that $\chi(\HH{n}{3})$ is not bounded above. We extend this result to $\chi(\hh{})$ for any fixed $r \ge 2$.

Since there exists a homomorphism from $\hh{}$ to $K_{n:r}$, we have that $\chi(\hh{}) \le \chi(K_{n:r})$. Since the chromatic number of the Kneser graph is known, we have the following theorem.

\begin{theorem}
For $n \ge 2r$, $\chi(\hh{}) \le n-2r+2$.
\end{theorem}

\proof
For $n \ge 2r$, $\chi(K_{n:r}) = n-2r+2$.\qed

The next bound we give shows that the chromatic number of \hh{} cannot behave too erratically with respect to $n$.

\begin{theorem}\label{recursivebd}
\[\chi(\HH{n-1}{r}) \le \chi(\hh{}) \le \chi(\HH{n-1}{r}) + 1.\]
\end{theorem}

\proof
The first inequality is trivial. The second inequality can be seen by examining the diagram of the 3-cell partition given above. We can color the $\HH{n-1}{r}$ subgraph of \hh{} with the colors $\{1,\ldots,\chi(\HH{n-1}{r})\}$, then we can color all of the vertices with $n$ in their tail (the middle cell) with $\chi(\HH{n-1}{r})+1$, then color all the vertices with $n$ as their head (the right cell) with color 1.\qed

So when we increase $n$ by one, the chromatic number either stays the same, or increases by one. Note that this recursive bound actually implies the bound we get from the homomorphism to $K_{n:r}$, since \HH{2r}{r} is bipartite and in this case $n-2r+2 = 2$. Also note that $\HH{2r+1}{r}$ has an odd cycle and thus chromatic number three. Therefore the $n-2r+2$ bound is always met for $n = 2r-1, 2r, \text{ and } 2r+1$ ($n=2r-1$ is an empty graph). We would hope that this bound is met with equality for all $n \ge 2r-1$ as it is with the Kneser graphs, however after computing some small examples we see that this is not the case.

\begin{table}[h]
\begin{center}
	\begin{tabular}{|c|c|c|c|}
    \hline $r$ & $n$ & $\chi(\hh{})$ & $n-2r+2$ \\
    \hline
		2 & 4 & 2 & 2 \\
		2 & 5 & 3 & 3 \\
		2 & 6 & 4 & 4 \\
		2 & 7 & 4 & 5 \\
		\hline
		3 & 6 & 2 & 2 \\
		3 & 7 & 3 & 3 \\
		\hline
	\end{tabular}
	\caption{Chromatic numbers of $\hh{}$.}
\end{center}
\end{table}

For $r = 2$, we see that the chromatic number does not go up by one as $n$ changes from 6 to 7. This of course means that $\chi(\HH{n}{2}) \le n-2r+1$ for $n \ge 7$ and thus the $n-2r+2$ bound cannot be met with equality for these values of $n$ and $r$.

So far the bounds we have given have been mostly upper bounds, except for the trivial lower bound given in Theorem~\ref{recursivebd}, thus it is still an open question as to whether $\chi(\hh{})$ is bounded in terms of $n$. We proceed to answer this question.

Knowing that the chromatic number of $\hh{}$ does not increase with every increase of $n$, it is natural to ask whether or not it is bounded by some finite value which depends only on $r$. Recall the result from \cite{amote} which states that a cubic graph admits a homomorphism to \HH{22}{3} if and only if it is triangle-free. A trivial extension of this result is that $r$-regular graphs admit a homomorphism into $\HH{n'}{r}$, where $n' = r\frac{(r-1)^3 - 1}{r-2} + 1$, if and only if they are triangle-free. Since triangle-free regular graphs can have arbitrarily large chromatic number, this means that we cannot hope to bound $\chi(\hh{})$ for all values of $r$ and $n$, but it may be possible to bound it for fixed $r$. However, this turns out to not be the case.

We begin with a lemma that relates the chromatic numbers of H\"{a}ggkvist--Hell graphs with different tail sizes.

\begin{lemma}\label{tailchi}
Let $r \ge 2$. For any $n$, define $n^* = \chi(\hh{})+n$. Then,
\[\chi(\hh{}) \le \chi(\HH{n^*}{r+1}).\]
\end{lemma}

\proof
We prove this by showing that $\hh$ is isomorphic to a subgraph of $\HH{n^*}{r+1}$. Let
\[f: V(\hh{}) \rightarrow \{n+1,\ldots,n+\chi(\hh)=n^*\}\]
be a proper coloring of $\hh{}$. Consider the map
\[g: \hh{} \rightarrow \HH{n^*}{r+1}\]
given by
\[g(h_u,T_u) = (h_u,T_u \cup \{f(u)\}).\]
It is easy to see that this is an injective homomorphism which proves the lemma.\qed

This result immediately allows us to prove this next vital lemma.

\begin{lemma}\label{unbddunbdd}
For a fixed $r \ge 2$, if $\chi(\hh{})$ is unbounded with respect to $n$, then $\chi(\HH{n}{r'})$ is unbounded with respect to $n$ for all $r' \ge r$.
\end{lemma}

\proof
It will suffice to show that it holds for $r' = r+1$. We will prove the contrapositive. Suppose that $\chi(\HH{n}{r+1})$ is bounded by $M_{r+1}$ with respect to $n$. Then by Lemma~\ref{tailchi}, we have that
\[\chi(\hh{}) \le \chi(\HH{\chi(\hh{})+n}{r+1}) \le M_{r+1}\]
for all $n$. Therefore $\chi(\hh{})$ is bounded with respect to $n$.\qed

All that is left is to prove that $\chi(\HH{n}{2})$ is not bounded with respect to $n$. This we proceed to do.

Define a graph $S_n$ as follows: the vertices of $S_n$ are all 3-element subsets of $[n]$, and two such subsets, say $\{x_1,x_2,x_3\}$ with $x_1<x_2<x_3$, and $\{y_1,y_2,y_3\}$ with $y_1<y_2<y_3$, are adjacent if $x_2 = y_1$ and $x_3 = y_2$. Note that $S_n$ is a directed graph, but we will also use $S_n$ to refer to its underlying undirected graph. It follows from Ramsey's theorem for the partition of triples \cite{ramseycolor} that the chromatic number of $S_n$ may be arbitrarily large if $n$ is large. This fact is key to the proof of the analogous result for \HH{n}{2}:

\begin{lemma}\label{S_n}
$\chi(S_n) \le \chi(\HH{n}{2}) \le \chi(S_n)+2$.
\end{lemma}

\proof
For a vertex $\{x_1,x_2,x_3\} \in V(S_n)$ with $x_1<x_2<x_3$, we let
\[f(x_1,x_2,x_3) = (x_2,\{x_1,x_3\}) \in V(\HH{n}{2}).\]
Clearly, this is injective. We claim that it is an injective homomorphism. Suppose that $x = \{x_1,x_2,x_3\}$ with $x_1<x_2<x_3$ and $y = \{y_1,y_2,y_3\}$ with $y_1<y_2<y_3$ are adjacent in $S_n$. Then WLOG $x_2 = y_1$ and $x_3 = y_2$. Now $f(x) = (x_2,\{x_1,x_3\})$ and $f(y) = (y_2,\{y_1,y_3\})$, and $x_2 \in \{y_1,y_3\}$, $y_2 \in \{x_1,x_3\}$ and $\{x_1,x_3\} \cap \{y_1,y_3\} = \varnothing$. Therefore $f(x)$ and $f(y)$ are adjacent. This proves that $S_n$ is isomorphic to a subgraph of \HH{n}{2} which implies the first inequality.

Now we will show that it is isomorphic to an induced subgraph of \HH{n}{2}. Suppose that $(x_2,\{x_1,x_3\})$ with $x_1<x_2<x_3$ and $y = (y_2,\{y_1,y_3\})$ with $y_1<y_2<y_3$ are adjacent in \HH{n}{2}. Either $x_2 = y_1$ or $x_2 = y_3$. Suppose that $x_2 = y_1$. Then $y_2 = x_3$ since $y_2>y_1=x_2$, so we have that $x_2 = y_1$ and $x_3 = y_2$. Therefore $\{x_1,x_2,x_3\}$ and $\{y_1,y_2,y_3\}$ are adjacent in $S_n$. Similarly, if $x_2 = y_3$, we deduce that $y_2 = x_1$, and thus $\{x_1,x_2,x_3\}$ and $\{y_1,y_2,y_3\}$ are again adjacent in $S_n$. Therefore, $S_n$ is isomorphic to the subgraph of \HH{n}{2} induced by the vertices whose heads are in between the two elements in their tails. The remaining vertices of \HH{n}{2} are the vertices whose heads are greater/less than both elements in their tails. This is simply the disjoint union of two independent sets of the recursive type, thus we can color the rest of \HH{n}{2} with two colors, and therefore we can color \HH{n}{2} with $\chi(S_n)+2$ colors.\qed

This immediately gives us that for any fixed $r \ge 2$ the chromatic number of \hh{} is unbounded with respect to $n$. Combining this with Lemma~\ref{recursivebd} we obtain the following:

\begin{theorem}\label{allunbddstrong}
For any $r\ge 2$, for any positive integer $k$, there exists an integer $n$ such that $\chi(\HH{n}{r}) = k$.\qed
\end{theorem}

It is important to note that in 2000, Gallucio, Hell, and Ne\u{s}et\u{r}il showed in \cite{color3} that the chromatic number of \HH{n}{3} is unbounded with respect to $n$. Their proof used a similar idea to our proof of Lemma~\ref{tailchi} in order to show that $S_n$ was a subgraph of \HH{n'}{3} for some $n' \ge n$. They also proved that $\chi(\HH{n}{3}) \ge 4$ for $n \ge 16$, however our computation of $\chi(\HH{6}{2}) = 4$ combined with Lemma~\ref{tailchi} proves the same but for $n \ge 10$.

So we have shown that $\chi(\hh{}) \le n-2r+2$, but that this bound is not always met with equality, that increasing $n$ by one increases the chromatic number by at most one, and that for any fixed $r \ge 2$ the chromatic number of \hh{} is not bounded with respect to $n$.

\section{Fractional Chromatic Number}
\label{sec:frac}

Let $\mathcal{I}(X)$ denote the set of all independent sets of the graph $X$, and $\mathcal{I}(X,x)$ the set of all independent sets of $X$ that contain $x$. A \emph{fractional coloring} of $X$ is defined to be a non-negative real-valued function $f$ on $\mathcal{I}(X)$ such that for any vertex $x$ of $X$,
\[\sum_{S\in \mathcal{I}(X,x)} f(S) \ge 1.\]
The \emph{weight} of a fractional coloring is the sum of all of its values, and the \emph{fractional chromatic number} of the graph $X$ is defined to be the minimum possible weight of a fractional coloring of $X$, and is denoted by $\chi^*(X)$.

As with the chromatic number, if $X$ and $Y$ are graphs and $X\rightarrow Y$, then $\chi^*(X) \le \chi^*(Y)$. Recall that coloring a graph can be seen as a homomorphism to a complete graph, and so the chromatic number of a graph $X$ can be defined as the minimum $n$ such that $X\rightarrow K_n$. For fractional chromatic number Kneser graphs play an analogous role to complete graphs and we have the following theorem:

\begin{theorem}\label{kneserchi*}
For any graph $X$ we have
\[\chi^*(X) = \min\{n/r : X \rightarrow K_{n:r}\}.\qed\]
\end{theorem}

For a vertex transitive graph $X$, the fractional chromatic number can be easily determined if the independence number is known:
\begin{theorem}\label{chi*trans}
If $X$ is a vertex transitive graph, then $\chi^*(X) = \frac{|V(X)|}{\alpha(X)}$.\qed
\end{theorem}
Using this and our lower bound on $\alpha(\hh{})$, we are able to give an upper bound on $\chi^*(\hh{})$. Combining this upper bound and Theorem~\ref{kneserchi*}, we deduce that for $n > r^2 + 1$ there exists a homomorphism from \hh{} to some Kneser graph which is not $K_{n:r}$. We then show that we are able to construct this homomorphism using the independent set from Theorem~\ref{bestindybd}.

We do not give the details of the above results concerning fractional chromatic number, but if the reader is interested Godsil and Royle's \textit{Algebraic Graph Theory} \cite{agt_godsil} is a good reference.

Recall from Theorem~\ref{bestindybd} that $\alpha(\hh{}) \ge \frac{r}{n}|V(\hh{})|$ for $n \ge 2r$, and $\alpha(\hh{}) \ge \frac{1}{r+1}|V(\hh{})| + \frac{r-1}{r+1}\binom{r^2}{r}$ for $n \ge r^2$. Combining this with \ref{chi*trans} gives us the following theorem:

\begin{theorem}\label{chi*upbd2}
For $n \ge 2r$,
\[\chi^*(\hh{}) \le \frac{n}{r}\]
For $n \ge r^2$,
\[\chi^*(\hh{}) \le (r+1)\left(1-\frac{(r-1)\binom{r^2}{r}}{(r+1)\binom{n}{r+1} + (r-1)\binom{r^2}{r}}\right) < r + 1\]
\end{theorem}

\proof
Arithmetic.\qed

Note that the second bound given above is strictly less than the first bound for $n > r^2 + 1$. This is because $\frac{1}{r+1}|V(\hh{})| + \frac{r-1}{r+1}\binom{r^2}{r} > \frac{r}{n}|V(\hh{})|$ for $n > r^2 + 1$. So we see that $\chi^*(\hh{}) < \frac{n}{r}$ for these values of $n$. Therefore, by Theorem~\ref{kneserchi*}, there must be some Kneser graph $K_{n':r'}$ with $\frac{n'}{r'} < \frac{n}{r}$ such that $\hh{} \rightarrow K_{n':r'}$ for $n \ge r^2 + 2$. Furthermore, we must have that $\frac{n'}{r'} < r+1$. So unlike the chromatic number, the fractional chromatic number of \hh{} is bounded for fixed $r$. This differs from the Kneser graphs which have fractional chromatic number $\chi^*(K_{n:r}) = n/r$ which is of course unbounded for fixed $r$.

In fact we need not speculate about the homomorphism whose existence is implied by the above, since we are able to construct it using the same result about independent sets of \hh{} that we used to prove Theorem~\ref{chi*upbd2}.

\begin{theorem}\label{frachom}
Let $X$ be a graph with independent set $S$, and let $G$ be a subgroup of $\aut(X)$ which acts transitively on $V(X)$. Then $X \rightarrow K_{n':r'}$, where $n' = |G|$ and $r' = \frac{|G||S|}{|V(X)|}$.
\end{theorem}

\proof
This proof is basically a special case of the proof of Theorem~\ref{kneserchi*} from \cite{agt_godsil}. To prove the theorem we will give a homomorphism to a graph which is isomorphic to $K_{n':r'}$. Let $Y$ be the graph whose vertices are the $r'$-subsets of the elements of $G$, and adjacency is disjointedness. Clearly this is isomorphic to $K_{n':r'}$. Now consider the map $\varphi: V(X) \rightarrow V(Y)$ defined as follows:
\[\varphi(x) = \{g \in G:x \in g(S)\}\]
Before we show that this is a graph homomorphism, we must first show that it actually does map to $V(Y)$. To show this we simply need to show that the set in the equation above has size $r'$ for every vertex of $X$. Since $G$ acts vertex transitively on $X$, it is clear that these sets will all have the same size, say $d$. To determine the size of these sets we can consider the hypergraph whose vertices are $V(X)$ and whose edges are $\{g(S) : g \in G\}$. The total degree of the vertices of this hypergraph is $d|V(X)|$, whereas the total degree of the edges is $|G||S|$. Since these two values must be equal, we have that $d = \frac{|G||S|}{|V(X)|} =r'$. So $\varphi$ is map to $V(Y)$. Now suppose that $x_1$ and $x_2$ are adjacent in $X$, then they are not both contained in any independent set of $X$ and in particular they are not contained in any image of $S$ under an element of $G$. Therefore there does not exist $g \in G$ such that $g \in \varphi(x_1)$ and $g \in \varphi(x_2)$, and so $\varphi(x_1) \cap \varphi(x_2) = \varnothing$ which means that $\varphi(x_1)$ is adjacent to $\varphi(x_2)$ in $Y$. Thus $\varphi$ is a graph homomorphism.\qed

Applying this theorem to \hh{} we obtain the following corollary:

\begin{corollary}
For $n \ge r^2$, $\hh{} \rightarrow K_{n':r'}$ where $n' = n!$ and
\begin{align*}
r' &= \left[\binom{n}{r+1} + \frac{r-1}{r+1}\binom{r^2}{r}\right](n-r-1)!r! \\
&= \frac{1}{r+1}\left[n! + (r-1)(r^2\cdot \ldots \cdot (r^2 - r + 1))(n-r-1)!\right]
\end{align*}
\end{corollary}

\proof
Use the independent set given in Theorem~\ref{bestindybd} and Theorem~\ref{frachom}.\qed

Using the independence numbers we computed for small $n$ and $r$, we can compute the fractional chromatic numbers for the same values:

\begin{center}
	\begin{tabular}{|c|c|c|}
    \hline $r$ & $n$ & $\chi^*(\hh{})$  \\
    \hline
		2 & 4 & 2 \\
		2 & 5 & 5/2 \\
		2 & 6 & 30/11 \\
		2 & 7 & 105/37 \\
		2 & 8 & 84/29 \\
		\hline
		3 & 6 & 2 \\
		3 & 7 & 7/3 \\
		3 & 8 & 8/3 \\
		\hline
	
	\end{tabular} \\
	\textsc{Table 3.} Fractional chromatic numbers of \HH{n}{2} and \HH{n}{3}.
\end{center}

\section{Automorphisms}
\label{sec:aut}

As we noted in \sref{basic}, the permutations on $n$ elements act as automorphisms of \hh{}, i.e. $S_n \subseteq \aut(\hh{})$. The aim of this section is to show that these are the only automorphisms of \hh{} for $n \ge 2r+1$. For $n \le 2r$, it is fairly easy to see that this is not the case, and so it is assumed that $n \ge 2r+1$ throughout this section.

The main difficulty in proving the above result is to show that any automorphism of \hh{} must preserve the property of two vertices having the same tail. To do this it suffices to show that it is possible to tell when a pair of vertices has the same tail. With this done, we can continue our trend of standing on the shoulders of the analogous result for Kneser graphs to finish the proof.

At this point it is convenient to define some terms. If two distinct vertices of \hh{} have the same tail, we refer to them as a \emph{tail-type pair}. If they have the same head and their tails have exactly $r-1$ elements in common, then we refer to them as a \emph{head-type pair}. The next four lemmas show that it is always possible to distinguish a tail-type pair from a pair of vertices which is not of tail-type.

In the next four lemmas, $x$ and $y$ will always be a tail-type pair and have the forms $(h_x,T)$ and $(h_y,T)$ respectively. Similarly, $z$ and $w$ will always be a head-type pair and have the forms $(h,C \cup \{t_z\})$ and $(h,C \cup \{t_w\})$ respectively. Here $h_x \ne h_y$ and $t_z \ne t_w$.

\begin{lemma}
\label{dominate}
Let $u,v \in V(\hh{})$ be distinct. If $u$ and $v$ are a tail-type pair, then there is no vertex of \hh{} which is adjacent to all of the common neighbors of $u$ and $v$. Furthermore, if $u$ and $v$ are neither a tail-type pair nor a head-type pair, then there is some other vertex of \hh{} which is adjacent to all of their common neighbors.
\end{lemma}

\proof
Let $x$ and $y$ be defined as above. Let $N_{xy}$ be the set of common neighbors of $x$ and $y$. Note that $N_{xy}$ is always nonempty. We claim that the only vertices adjacent to every vertex in $N_{xy}$ are $x$ and $y$.

Let $s = (h_s,T_s)$ be a vertex of \hh{} that is not $x$ or $y$. If $T_s \ne T$, then let $h' \in T\setminus T_s$ and $h_x,h_y \in T' \subseteq [n]\setminus T$. Then $(h',T') \in N_{xy}$, but $s$ is not adjacent to $(h',T')$. Now if $T_s = T$, then $h_x \ne h_s \ne h_y$, so let $h' \in T$ and $h_x,h_y \in T' \subseteq [n] \setminus (T \cup \{h_s\})$. Then again we have that $(h',T') \in N_{xy}$ but $s$ is not adjacent to $(h',T')$. So we have proved the claim.

It is both necessary and straightforward to show that the choices for $h'$ and $T'$ above are always possible whenever $n \ge 2r+1$.

Now suppose that $u,v \in V(\hh{})$ are neither a tail-type pair nor a head-type pair. We will show that there is a vertex $s$, distinct from $u$ and $v$, which is adjacent to every common neighbor of $u$ and $v$. There are two main cases: $h_u \ne h_v$ and $h_u = h_v$. We begin with the case of unequal heads.

Let $s = (h_u,T_v)$. Now let $(h',T')$ be a common neighbor of $u$ and $v$. This implies that $h' \in T_v$, $h_u \in T'$, and $T' \cap T_v = \varnothing$, and so $(h',T')$ must also be adjacent to $s$. Therefore $s$ is adjacent to every common neighbor of $u$ and $v$ and so we are done.

Now consider the case where $h_u = h = h_v$. Let $D = T_u \cap T_v$. If $D = \varnothing$, then there are no common neighbors of $u$ and $v$ since there are no choices for the head of a such a vertex. So $D$ can be assumed nonempty, and since $u$ and $v$ are not a head-type pair, we have that $1 \le |D| \le r-2$. Let $t_u \in T_u \setminus T_v$ and $t_v \in T_v \setminus T_u$. Let $s = (h,T_s)$ where $T_s = T_u \cup \{t_v\} \setminus \{t_u\} \supseteq D$, and note that $T_s$ is not $T_v$ as it would be in the case where $|D| = r-1$. If $(h',T')$ is a common neighbor of $u$ and $v$, then $h' \in D$, $h \in T'$, and $T' \cap (T_u \cup T_v) = \varnothing$. So we see that $s$ is also adjacent to $(h',T')$ and so we are done.\qed

So we have only left to show that we can distinguish between a tail-type pair and a head-type pair. The next three lemmas take care of this case.

\begin{lemma}
\label{common}
Let $x,y \in V(\hh{})$ be a tail-type pair, and let $z,w \in V(\hh{})$ be a head-type pair. The number of common neighbors of $x$ and $y$ is equal to the number of common neighbors of $z$ and $w$ if and only if $n = 3r$.
\end{lemma}

\proof
Recall the set $N_{xy}$ of common neighbors of $x$ and $y$. It is straightforward to see that
\[N_{xy} = \left\{(h',T') \in V(\hh{}) : h' \in T, \ h_x,h_y \in T', \ \& \ T' \cap T = \varnothing\right\}.\]
Considering a possible vertex $(h',T')$ in $N_{xy}$, we see that
\[T' \setminus \{h_x,h_y\} \subseteq [n]\setminus (T \cup \{h_x,h_y\})\]
which means that there are $\binom{n-r-2}{r-2}$ choices for the tail, and there are $|T| = r$ choices for its head. So there are $r\binom{n-r-2}{r-2}$ common neighbors of $x$ and $y$.

Now consider the set $N_{zw}$ of common neighbors of $z$ and $w$. Again it is straightforward to see that
\[N_{zw} = \left\{(h',T') \in V(\hh{}) : h' \in C, \ h \in T', \ \& \ T' \cap (C \cup \{t_z,t_w\}) = \varnothing \right\}.\]
Again considering a possible vertex $(h',T')$ in $N_{zw}$, we see that
\[T'\setminus \{h\} \subseteq [n]\setminus (C \cup \{h,t_z,t_w\})\]
giving us $\binom{n-r-2}{r-1}$ choices for the tail, and $|C|=r-1$ choices for the head. So there are $(r-1)\binom{n-r-2}{r-1}$ common neighbors of $z$ and $w$.

So the number of common neighbors of $x$ and $y$ is equal to the number of common neighbors of $z$ and $w$ if and only if $r\binom{n-r-2}{r-2}=(r-1)\binom{n-r-2}{r-1}$. Simple arithmetic shows that this is equivalent to $n = 3r$.\qed

Lemmas~\ref{dominate} and \ref{common} prove that we are able to tell when a pair of vertices has the same tail in every case except for $n=3r$. For this case we need the next two lemmas.

\begin{lemma}
\label{uncommon}
Let $r \ge 3$ and $n \ge 3r-1$. Let $x,y \in V(\hh{})$ be a tail-type pair, and let $z,w \in V(\hh{})$ be a head-type pair. If a vertex of \hh{} is adjacent to a neighbor of $x$ and a neighbor of $y$, then it is adjacent to some common neighbor of $x$ and $y$. The same is not true for $z$ and $w$.
\end{lemma}

\proof
Let $v = (h_v,T_v)$ be a vertex that is adjacent to both a neighbor of $x$ and a neighbor of $y$. Also, let $T_{x'}$ and $T_{y'}$ be the tails of these neighbors respectively. Since any neighbor of $x$ or $y$ must have an element of $T$ as its head, we can conclude that $v$ has an element, say $t$, of $T$ in its tail. Furthermore, since any neighbor of $x$ must have $h_x$ in its tail, $v$ must not contain $h_x$ in its tail, and similarly for $h_y$. Now $h_v \in T_{x'} \cap T_{y'}$ and $|T_{x'} \cup T_{y'} \setminus \{h_x,h_y\}| \ge r-2$, so let $T'$ be a subset of this with size $r-2$ which contains $h_v$. Then $(t,T'\cup \{h_x,h_y\})$ is adjacent to $x$, $y$, and $v$. Note that this does not work for $r=2$ since $T'$ would have to both contain $h_v$ and have size 0 in this case.

Now consider the vertices $u_z = (t_z,\{h\} \cup T')$ and $u_w = (t_w,\{h\} \cup T')$ where $T' \subseteq [n]\setminus (C \cup \{h,t_z,t_w\})$ and $|T'| = r-1$. Note that these are neighbors of $z$ and $w$ respectively, but neither is a neighbor of both. Now let $T_v \subseteq [n]\setminus (C \cup T' \cup \{h,t_z,t_w\})$ such that $|T_v| = r-2$. Notice that
\[\big|[n]\setminus (C \cup T' \cup \{h,t_z,t_w\})\big| = n - 2r - 1 \ge r-2,\]
and so our choice of $T_v$ is possible. Letting $v = (h,T_v \cup \{t_z,t_w\})$, we see that $v$ is adjacent to both $u_z$ and $u_w$. However, since every common neighbor of $z$ and $w$ has an element of $C$ as its head, $v$ is adjacent to none of them and so we are done.\qed

We have now taken care of every case except for when $r=2$ and $n=6$. This final case is resolved with the following lemma.

\begin{lemma}
\label{uncommon2}
Let $r=2$ and $n \ge 6$. Let $x,y \in V(\hh{})$ be a tail-type pair, and let $z,w \in V(\hh{})$ be a head-type pair. If a vertex is adjacent to a neighbor of $x$ which is not a neighbor of $y$, and a neighbor of $y$ which is not a neighbor of $x$, then it is not adjacent to any common neighbor of $x$ and $y$. The same is not true for $z$ and $w$.
\end{lemma}

\proof
Let $N_x$ be the set of neighbors of $x$ which are not also neighbors of $y$, and define $N_y$ similarly. The vertices in $N_x$ have an element of $T$ for their head, and tails disjoint from $T$. So no vertex in $N_x$ has $h_y$ in its tail, otherwise it would also be a neighbor of $y$. Similarly, no vertex in $N_y$ has $h_x$ in its tail. Suppose $v$ is adjacent to a vertex in $N_x$ and a vertex in $N_y$. Then $v$ cannot have $h_x$ as its head, as no vertex of $N_y$ has $h_x$ in its tail, and similarly $v$ cannot have $h_y$ as its head. But since $r=2$, the tails of all of the common neighbors of $x$ and $y$ are exactly $\{h_x,h_y\}$, and so $v$ is adjacent to none of them. So if a vertex is adjacent to some vertex in $N_x$ and some vertex in $N_y$, then it is not adjacent to any common neighbor of $x$ and $y$. All that is left is to show that this does not hold for $z$ and $w$.

In this case $|C| = 1$, so let $c$ be the only element in $C$. Let $t_1,t_2 \in [n]\setminus \{h,c,t_z,t_w\}$, which is possible since $n \ge 6$. The vertices $(c,\{h,t_w\})$, $(c,\{h,t_z\})$, and $(c,\{h,t_1\})$ are adjacent to only $z$, only $w$, and both respectively. The vertex $(h,\{c,t_2\})$ is adjacent to all three of these vertices and so we are done.\qed

So we have finally shown that a pair of vertices with the same tail is distinguishable from a pair of vertices with different tails. This immediately gives us the next theorem, which will allow us to determine the automorphism group of the H\"{a}ggkvist--Hell graphs.

\begin{theorem}
\label{tailaut}
Let $\varphi \in \aut(\hh{})$, and $x,y \in V(\hh{})$. If $x$ and $y$ have the same tail, then so do $\varphi(x)$ and $\varphi(y)$.
\end{theorem}

\proof
Lemmas \ref{dominate} through \ref{uncommon2}.\qed

What this theorem implies is that an automorphism of \hh{} gives rise to a bijection on the set of all possible tails. This can be thought of as a bijection on the vertices of the Kneser graph $K_{n:r}$, and can be shown to be an automorphism of that graph.

\begin{lemma}
\label{fixtrivial}
Let $\varphi \in \aut(\hh{})$. If $\varphi$ fixes the tail of every vertex in \hh{}, then $\varphi$ is identity.
\end{lemma}

\proof
Suppose not. Then there exists a vertex $v = (h_v,T)$ such that $\varphi(v) = u = (h_u,T)$ where $h_u \ne h_v$. Let $T'$ be an $r$-subset of $[n]\setminus(T \cup \{h_u\})$ containing $h_v$, and let $t \in T$. Then $v$ is adjacent to $(t,T')$, but $u$ is not adjacent to any vertex with tail $T'$, since $h_u \notin T'$, and thus it is not adjacent to $\varphi(t,T')$, since $\varphi$ fixes tails. This is a contradiction of the definition of automorphism and so we are done.\qed

Note that our choice of $T'$ above is possible since $n \ge 2r+1$. For $n = 2r$ the above lemma does not hold.

In order to prove that the automorphism group of \hh{} is isomorphic to $S_n$, we need to enlist the analogous result for the Kneser graph $K_{n:r}$. We give this now without proof, but if the reader is interested the proof relies on Erd\"{o}s-Ko-Rado and can be found in \cite{agt_godsil}.

\begin{theorem}
\label{kneserautos}
For $n \ge 2r+1$, $\aut(K_{n:r}) \cong S_n$.\qed
\end{theorem}

We are now able to prove the main result of this section.

\begin{theorem}
\label{automorphisms}
For $n \ge 2r+1$, $\aut(\hh{}) \cong S_n$.
\end{theorem}

\proof
Let $\varphi \in \aut(\hh{})$. Define $\varphi^*$ to be the bijection on the $r$-subsets of $[n]$ given by the following: if $\varphi(h_1,T_1) = (h_2,T_2)$, then $\varphi^*(T_1) = T_2$. By Theorem~\ref{tailaut} $\varphi^*$ is well defined, and it is a bijection since otherwise $\varphi$ would not be. We will show that $\varphi^*$ is a automorphism of $K_{n:r}$.

Let $S$ and $T$ be two adjacent vertices of $K_{n:r}$, i.e.~two disjoint $r$-subsets of $[n]$. Pick $s \in S$ and $t \in T$. Then $(t,S)$ is adjacent to $(s,T)$ in \hh{} and so $\varphi(t,S) = (t',S')$ is adjacent to $\varphi(s,T) = (s',T')$. But this means that $S' = \varphi^*(S)$ is disjoint from $T' = \varphi^*(T)$, and thus they are adjacent in $K_{n:r}$. Therefore $\varphi^* \in \aut(K_{n:r}) \cong S_n$. Let $\sigma$ be the inverse of $\varphi^*$ in $S_n$. Then $\sigma$ is an automorphism of \hh{}, and therefore $\sigma \circ \varphi$ is an automorphism of \hh{} and it clearly fixes the tail of every vertex and so by Lemma~\ref{fixtrivial} it is the identity. Therefore $\varphi = \sigma^{-1} \in S_n$ and we are done.\qed

%
%

\section{Discussion and Open Questions}

This paper has shown that the relationship between \hh{} and $K_{n:r}$ is both strong and useful. This relationship was key in determining the diameter and odd girth of \hh{}, and these parameters are equal to those of $K_{n:r}$ for $n < \frac{7}{3}r - \frac{1}{3}$ and $n < 3r$ respectively. Similarly, the bound given for $\alpha(\hh{})$ relies on the value of $\alpha(K_{n:r})$, and is in fact equal for $n \le r^2 + 1$. Though we did not prove it, we believe that this bound is met with equality for these values of $n$. Again with the chromatic number we were able to use the homomorphism from \hh{} to $K_{n:r}$ to give an upper bound. Furthermore, we were able to show that, like the Kneser graphs, the chromatic number of \hh{} is unbounded for fixed $r$.

However, we have also seen throughout this paper that for fixed $r$, as $n$ increases the similarities between \hh{} and $K_{n:r}$ tend to fade away and the gaps between parameters of these graphs grow larger. The diameter of Kneser graphs eventually becomes two as $n$ increases, whereas the diameter of H\"{a}ggkvist--Hell graphs is always at least four. Similarly for odd girth, which is three for $K_{n:r}$ with $n \ge 3r$ but is never less than five for \hh{}. This reminds us that though H\"{a}ggkvist--Hell graphs are triangle-free, and thus have clique number two, the Kneser graphs have $\omega(K_{n:r}) = \left\lfloor n/r \right\rfloor$. This difference is possibly what is being expressed by the fact that for fixed $r$, H\"{a}ggkvist--Hell graphs have bounded fractional chromatic number but Kneser graphs do not ($\chi^*(X) \ge \omega(X)$ for all graphs $X$). We note that $|V(\hh{})| = (n-r)|K_{n:r}|$ and $|E(\hh{})| = r^2|E(K_{n:r})|$, which means that as $n$ increases, the H\"{a}ggkvist--Hell graphs become sparser and sparser relative to the Kneser graphs. It makes sense then that the diameter, odd girth, and independence number of \hh{} grow larger in comparison to these parameters for $K_{n:r}$ as $n$ increases, since these parameters are typically larger for sparser graphs. And similarly for the relationship between the fractional chromatic numbers of \hh{} and $K_{n:r}$.

In this paper we have studied only a handful of graph parameters which are typically of interest to graph theorists, so of course there are many more questions to be asked about these graphs. Here are some that we find worth pondering.

Any improvement on bounds for the independence number, chromatic number, or fractional chromatic number, would be of interest, especially if one were able determine the exact value of any of them. The independence number seems likely to be the easiest candidate for the latter, though of course this immediately gives us the fractional chromatic number.

There are quite a few results regarding homomorphisms between different Kneser graphs; we would be curious to see which of these are able to be extended to H\"{a}ggkvist--Hell graphs. In particular, it is quite easy to see that \hh{} is an induced subgraph of $\HH{tn}{tr}$ for any positive integer $t$, as is analogously true for Kneser graphs. However, it is not so clear as to whether there is an analogous homomorphism to the one from $K_{n:r}$ to $K_{n-2:r-1}$ \cite{agt_godsil} for H\"{a}ggkvist--Hell graphs.

Though we did not go into the study of them, we are interested in whether or not H\"{a}ggkvist--Hell graphs are cores. Cores are graphs with no proper endomorphisms, and they are the minimal elements of the equivalence classes of homomorphic equivalence. Two graphs $X$ and $Y$ are homomorphically equivalent if $X \rightarrow Y$ and $Y \rightarrow X$. It is known that $K_{n:r}$ is a core for $n \ge 2r+1$. We would like to see if the same is true for H\"{a}ggkvist--Hell graphs. We suspect it is.

Finally, there are some quite natural generalizations of H\"{a}ggkvist--Hell graphs that are of potential interest. The most obvious generalization is to let the heads of vertices be of sizes other than one. In other words, the vertices are all ordered pairs $(\alpha, \beta)$ of subsets of $[n]$ where
\[|\alpha| = r_1, \ \ |\beta| = r_2, \ \ \alpha \cap \beta = \varnothing.\]
Two vertices $(\alpha,\beta)$ and $(\alpha',\beta')$ are adjacent if
\[\alpha \subseteq \beta', \ \ \alpha' \subseteq \beta, \ \ \beta \cap \beta' = \varnothing.\]
We can also consider the $q$-analogs of these graphs similarly to $q$-Kneser graphs. The $q$-Kneser graph, $qK_{n:r}$, has $r$-dimensional subspaces of an $n$-dimensional vector space over a finite field of order $q$ as vertices. Adjacency in $qK_{n:r}$ is having intersection equal to the 0-dimensional subspace. The $q$-analogs of \hh{} would be defined analogously.

\section*{Acknowledgements}
The author would like to thank his advisor, Chris Godsil, for the introduction to these graphs and guidance in the study of them. The author would also like to thank Krystal Guo for a helpful discussion that aided in determining the automorphism group of these graphs. Most of the results presented here are from the author's Master's thesis completed in the department of Combinatorics and Optimization at the University of Waterloo.

\renewcommand{\bibname}{References}

\bibliographystyle{plain}
\bibliography{haggkvist_hell}

\end{document}